\documentclass{amsart}

\usepackage{latexsym}
\usepackage{amsthm,amsfonts,amssymb,mathrsfs}
\usepackage{rotating}
\usepackage[leqno]{amsmath}
\usepackage{xspace}
\usepackage[all]{xy}
\usepackage{longtable}
\newcommand{\numberseries}{\mdseries}   %Fontseries used for numbering theorem

\newlength{\thmtopspace}                %Space above theorem
\newlength{\thmbotspace}                %Space below theorem
\newlength{\thmheadspace}               %Space between theorem caption and text
\newlength{\thmindent}                  %For indenting

\newtheoremstyle{bfupright head,slanted body}
                {\thmtopspace}{\thmbotspace}
                {\slshape}{\thmindent}{\bfseries}{.}{\thmheadspace}
                {{\numberseries \thmnumber{(#2) }}\thmnote{#3}}

\newtheoremstyle{bfupright head,upright body}
                {\thmtopspace}{\thmbotspace}
                {\upshape}{\thmindent}{\bfseries}{.}{\thmheadspace}
                {{\numberseries \thmnumber{(#2) }}\thmnote{#3}}

\newtheoremstyle{bfit head,upright body}
                {\thmtopspace}{\thmbotspace}
                {\upshape}{\thmindent}{\upshape}{.}{\thmheadspace}
                {{\numberseries\thmnumber{(#2) }}
                {\bfseries\itshape\thmnote{\negthickspace#3}}}

\newtheoremstyle{it head,upright body}
                {\thmtopspace}{\thmbotspace}
                {\upshape}{\thmindent}{\upshape}{.}{\thmheadspace}
                {{\numberseries\thmnumber{(#2) }}
                {\itshape\thmnote{\negthickspace#3}}}

\newtheoremstyle{fixed bf head,slanted body}
                {\thmtopspace}{\thmbotspace}{\slshape}
                {\thmindent}{\bfseries}{.}{\thmheadspace}
                {{\numberseries \thmnumber{(#2) }}\thmname{#1}\thmnote{ (#3)}}

\newtheoremstyle{fixed bf head,upright body}
                {\thmtopspace}{\thmbotspace}{\upshape}
                {\thmindent}{\bfseries}{.}{\thmheadspace}
                {{\numberseries \thmnumber{(#2) }}\thmname{#1}\thmnote{ (#3)}}

\newtheoremstyle{independent paragraph}
                {\thmtopspace}{\thmbotspace}
                {\upshape}{\thmindent}{\upshape}{}{0pt}
                {\thmnote{#3 }}

\newtheoremstyle{subparagraph}
                {\thmbotspace}{\thmbotspace}
                {\upshape}{\thmindent}{\upshape}{}{0pt}
                {\thmnote{#3 }}

\newtheoremstyle{notes}
                {\thmtopspace}{\thmbotspace}
                {\ttfamily}{\thmindent}{\ttfamily\small }{}{0pt}
                {\thmnote{#3 }}

\theoremstyle{bfupright head,slanted body}
\newtheorem{res}{}[section]             \newtheorem*{res*}{}

\theoremstyle{bfit head,upright body}
                 \newtheorem*{com*}{}

\theoremstyle{bfupright head,upright body}
\newtheorem{bfhpg}[res]{}               \newtheorem*{bfhpg*}{}

\theoremstyle{it head,upright body}
               \newtheorem*{ithpg*}{}

\theoremstyle{fixed bf head,slanted body}
\newtheorem{thm}[res]{Theorem}          \newtheorem*{thm*}{Theorem}
\newtheorem{prp}[res]{Proposition}      \newtheorem*{prp*}{Proposition}
\newtheorem{cor}[res]{Corollary}        \newtheorem*{cor*}{Corollary}
\newtheorem{lem}[res]{Lemma}            \newtheorem*{lem*}{Lemma}

\theoremstyle{fixed bf head,upright body}
       \newtheorem*{dfn*}{Definition}
\newtheorem{obs}[res]{Observation}      \newtheorem*{obs*}{Observation}
\newtheorem{rmk}[res]{Remark}           \newtheorem*{rmk*}{Remark}
\newtheorem{exa}[res]{Example}          \newtheorem*{exa*}{Example}
         \newtheorem*{exe*}{Exercise}
            \newtheorem{stp*}{Setup}

\theoremstyle{independent paragraph}
\newtheorem{ipg}{}

\theoremstyle{subparagraph}
\newtheorem{spg}{}

\theoremstyle{notes}

\newlength{\thmlistleft}        %leftmargin
\newlength{\thmlistright}       %rightmargin
\newlength{\thmlistpartopsep}   %partopsep
\newlength{\thmlisttopsep}      %topsep
\newlength{\thmlistparsep}      %parsep
\newlength{\thmlistitemsep}     %itemsep

\newcounter{eqc} 
  {\end{list}}%

\newcounter{prt}
\newenvironment{prt}{\begin{list}{\upshape (\alph{prt})}%
    {\usecounter{prt}%
      \setlength{\leftmargin}{\thmlistleft}%
      \setlength{\labelwidth}{\thmlistleft}%
      \setlength{\rightmargin}{\thmlistright}%
      \setlength{\partopsep}{\thmlistpartopsep}%
      \setlength{\topsep}{\thmlisttopsep}%
      \setlength{\parsep}{\thmlistparsep}%
      \setlength{\itemsep}{\thmlistitemsep}}}%
  {\end{list}}%

\newcommand{\prtlbl}[1]{{\upshape(#1)}}

\newcommand{\prevprt}[1]{{\upshape(\alph{prt}#1)}}

\newcounter{rqm}
\newenvironment{rqm}{\begin{list}{\upshape (\arabic{rqm})}%
    {\usecounter{rqm}%
      \setlength{\leftmargin}{\thmlistleft}%
      \setlength{\labelwidth}{\thmlistleft}%
      \setlength{\rightmargin}{\thmlistright}%
      \setlength{\partopsep}{\thmlistpartopsep}%
      \setlength{\topsep}{\thmlisttopsep}%
      \setlength{\parsep}{\thmlistparsep}%
      \setlength{\itemsep}{\thmlistitemsep}}}%
  {\end{list}}%

\newenvironment{eqlist}{\begin{list}{\upshape (\theequation)\hfill}%
    {\usecounter{equation}\setlength{\labelsep}{0pt}%
      \setlength{\leftmargin}{3.5em}%
      \setlength{\labelwidth}{3.5em}%
      \setlength{\rightmargin}{\thmlistright}%
      \setlength{\partopsep}{\thmlistpartopsep}%
      \setlength{\topsep}{\thmlisttopsep}%
      \setlength{\parsep}{\thmlistparsep}%
      \setlength{\itemsep}{\thmlistitemsep}}}%
  {\end{list}}%

\newenvironment{itemlist}{\nopagebreak \begin{list}{$\bullet$}%
    {\setlength{\leftmargin}{\thmlistleft}%
      \setlength{\labelwidth}{\thmlistleft}%
      \setlength{\rightmargin}{\thmlistright}%
      \setlength{\partopsep}{\thmlistpartopsep}%
      \setlength{\topsep}{\thmlisttopsep}%
      \setlength{\parsep}{\thmlistparsep}%
      \setlength{\itemsep}{\thmlistitemsep}}}%
  {\end{list}}%

\newenvironment{thmquote}[1][2.5em]{\begin{list}{}%
                        {\setlength{\leftmargin}{#1}\setlength{\rightmargin}{#1}%
                        \setlength{\partopsep}{0pt}%
                        \setlength{\topsep}{\thmbotspace}%
                        \setlength{\parsep}{0pt}%
                        \setlength{\itemsep}{0pt}}
                        \item[]}
                        {\end{list}}%

\newenvironment{prf}[1][Proof]{\begin{proof}[\bf #1]}{\end{proof}}

  \newcommand{\step}[1]{${#1}^\circ$}

  \newcommand{\proofoftag}[2][:]{(#2)#1}

\newcommand{\pgref}[1]{(\ref{#1})}
\newcommand{\pgpartref}[2]{(\ref{#1})\prtlbl{#2}}

\renewcommand{\eqref}[1]{\pgref{eq:#1}}

\newcommand{\corref}[2][Corollary~]{#1\pgref{cor:#2}}

\newcommand{\exaref}[2][Example~]{#1\pgref{exa:#2}}
\newcommand{\lemref}[2][Lemma~]{#1\pgref{lem:#2}}
\newcommand{\obsref}[2][Observation~]{#1\pgref{obs:#2}}
\newcommand{\prpref}[2][Proposition~]{#1\pgref{prp:#2}}
\newcommand{\resref}[1]{\pgref{res:#1}}

\newcommand{\thmref}[2][Theorem~]{#1\pgref{thm:#2}}

\newcommand{\corpartref}[3][Corollary~]{#1\pgpartref{cor:#2}{#3}}

\newcommand{\exapartref}[3][Example~]{#1\pgpartref{exa:#2}{#3}}
\newcommand{\lempartref}[3][Lemma~]{#1\pgpartref{lem:#2}{#3}}
\newcommand{\obspartref}[3][Observation~]{#1\pgpartref{obs:#2}{#3}}
\newcommand{\prppartref}[3][Proposition~]{#1\pgpartref{prp:#2}{#3}}
\newcommand{\respartref}[2]{\pgpartref{res:#1}{#2}}

\newcommand{\thmpartref}[3][Theorem~]{#1\pgpartref{thm:#2}{#3}}

\newcommand{\chpcite}[2][?]{\cite[ch.~#1]{#2}}
\newcommand{\corcite}[2][?]{\cite[cor.~#1]{#2}}
\newcommand{\lemcite}[2][?]{\cite[lem.~#1]{#2}}
\newcommand{\prpcite}[2][?]{\cite[prop.~#1]{#2}}
\newcommand{\rescite}[2][?]{\cite[#1]{#2}}
\newcommand{\seccite}[2][?]{\cite[sec.~#1]{#2}}
\newcommand{\thmcite}[2][?]{\cite[thm.~#1]{#2}}

\newcommand{\catbl}{\sqsubset}
\newcommand{\catbr}{\sqsupset}
\newcommand{\catb}{\sqsubset\mspace{-13mu}\sqsupset}

\newcommand{\Cat}[2]{{\sf{#2}}(#1)}
\newcommand{\Catsup}[3]{{\sf{#2}}^{\text{\upshape #3}}(#1)}
\newcommand{\Catsub}[3]{{\sf{#2}}_{#3}(#1)}
\newcommand{\Catsupsub}[4]{{\sf{#2}}^{\text{\upshape #3}}_{#4}(#1)}

\newcommand{\Catbl}[2]{\Catsub{#1}{#2}{\catbl}}
\newcommand{\Catbr}[2]{\Catsub{#1}{#2}{\catbr}}
\newcommand{\Catb}[2]{\Catsub{#1}{#2}{\catb}}

\newcommand{\Catbls}[3]{\Catsupsub{#1}{#2}{#3}{\catbl}}
\newcommand{\Catbrs}[3]{\Catsupsub{#1}{#2}{#3}{\catbr}}
\newcommand{\Catbs}[3]{\Catsupsub{#1}{#2}{#3}{\catb}}

\newcommand{\C}[1][R]{\Cat{#1}{C}}

\newcommand{\Cbr}[1][R]{\Catbr{#1}{C}}

\newcommand{\Cf}[1][R]{\Catsup{#1}{C}{\hspace{0.08em}f}}

\newcommand{\CfP}[1][R]{\Catsup{#1}{C}{fP}}
\newcommand{\CP}[1][R]{\Catsup{#1}{C}{P}}
\newcommand{\CF}[1][R]{\Catsup{#1}{C}{F}}
\newcommand{\CI}[1][R]{\Catsup{#1}{C}{I}}

\newcommand{\CfPbr}[1][R]{\Catbrs{#1}{C}{fP}}
\newcommand{\CPbr}[1][R]{\Catbrs{#1}{C}{P}}

\newcommand{\CIbl}[1][R]{\Catbls{#1}{C}{I}}

\newcommand{\CPb}[1][R]{\Catbs{#1}{C}{P}}
\newcommand{\CFb}[1][R]{\Catbs{#1}{C}{F}}
\newcommand{\CIb}[1][R]{\Catbs{#1}{C}{I}}

\renewcommand{\P}[1][R]{\Cat{#1}{P}}
\newcommand{\F}[1][R]{\Cat{#1}{F}}
\newcommand{\I}[1][R]{\Cat{#1}{I}}
\newcommand{\GP}[1][R]{\Cat{#1}{GP}}
\newcommand{\GI}[1][R]{\Cat{#1}{GI}}
\newcommand{\GF}[1][R]{\Cat{#1}{GF}}

\newcommand{\Pf}[1][R]{\Catsup{#1}{P}{\hspace{0.02em}f}}
\newcommand{\If}[1][R]{\Catsup{#1}{I}{\hspace{0.03em}f}}

\newcommand{\GPf}[1][R]{\Catsup{#1}{GP}{\hspace{0.02em}f}}

\newcommand{\GFf}[1][R]{\Catsup{#1}{GF}{\hspace{0.05em}f}}

\newcommand{\A}[1][R]{\Cat{#1}{A}}
\newcommand{\B}[1][R]{\Cat{#1}{B}}

\renewcommand{\H}[2][\no]{\operatorname{H}_{#1}(#2)}

\newcommand{\Tsr}[2]{#2\mspace{-2mu}\sideset{_{#1}}{}{\operatorname{\supset}}}

\newcommand{\Hom}[3][R]{\operatorname{Hom}_{#1}(#2,#3)}

\newcommand{\tp}[3][R]{#2\otimes_{#1}#3}

\newcommand{\Tor}[4][R]{\operatorname{Tor}^{#1}_{#2}(#3,#4)}

\newcommand{\bett}[3][R]{\beta^{#1}_{#2}(#3)}

\newcommand{\D}[1][R]{\Cat{#1}{D}}
\newcommand{\Dbl}[1][R]{\Catbl{#1}{D}}
\newcommand{\Dbr}[1][R]{\Catbr{#1}{D}}
\newcommand{\Db}[1][R]{\Catb{#1}{D}}

\newcommand{\Dfbr}[1][R]{\Catbrs{#1}{D}{f}}
\newcommand{\Dfb}[1][R]{\Catbs{#1}{D}{f}}

\newcommand{\DHom}[3][R]{\operatorname{\mathbf{R}Hom}_{#1}(#2,#3)}
\newcommand{\Dtp}[3][R]{#2\otimes_{#1}^{\mathbf{L}}#3}

\newcommand{\DG}[2][\mathfrak{a}]{\mbox{\ensuremath{\mathbf{R}\Gamma_{#1} #2}}}
\newcommand{\DL}[2][\mathfrak{a}]{\mbox{\ensuremath{\mathbf{L}\Lambda^{#1} #2}}}

\newcommand{\f}{\varphi}

\renewcommand{\l}{\ell}

\newcommand{\eq}{\simeq}
\newcommand{\is}{\cong}

\newcommand{\ZZ}{\mathbb{Z}}
\newcommand{\QQ}{\mathbb{Q}}

\newcommand{\m}{\mathfrak{m}}
\newcommand{\n}{\mathfrak{n}}
\newcommand{\p}{\mathfrak{p}}
\newcommand{\q}{\mathfrak{q}}

\newcommand{\Rmk}{(R,\m,k)}
\newcommand{\Snl}{(S,\n,l)}

\newcommand{\Shat}{\widehat{S}}

\newcommand{\onto}{\twoheadrightarrow}
\newcommand{\into}{\hookrightarrow}

\newcommand{\xra}{\xrightarrow}

\newcommand{\xre}{\xra{\;\eq\;}}

\newcommand{\mapdef}[4][\rightarrow]{\mbox{\ensuremath{#2\!: #3 #1 #4}}}
\newcommand{\dmapdef}[4][\lora]{#2:\; #3\:#1\:#4}

\newcommand{\poly}[2][k]{#1[#2]}
\newcommand{\pows}[2][k]{#1[\mspace{-2.3mu}[#2]\mspace{-2.3mu}]}

\newcommand{\supremum}[2]{\sup{\{\,#1\:|\:#2\}}}

\newcommand{\E}[2][R]{\operatorname{E}_{#1}(#2)}
\newcommand{\ERp}[1][1]{\E{R/\p}}

\newcommand{\Ker}[1]{\mbox{\ensuremath{\operatorname{Ker}#1}}}
\newcommand{\Coker}[1]{\mbox{\ensuremath{\operatorname{Coker}#1}}}

\newcommand{\tev}[1]{\omega_{#1}}
\newcommand{\hev}[1]{\theta_{#1}}

\newcommand{\dptR}{\operatorname{depth}R}

\newcommand{\fd}[2][R]{\operatorname{fd}_{#1}#2}
\newcommand{\id}[2][R]{\operatorname{id}_{#1}#2}
\newcommand{\pd}[2][R]{\operatorname{pd}_{#1}#2}

\newcommand{\Gfd}[2][R]{\operatorname{Gfd}_{#1}#2}
\newcommand{\Gid}[2][R]{\operatorname{Gid}_{#1}#2}
\newcommand{\Gpd}[2][R]{\operatorname{Gpd}_{#1}#2}

\newcommand{\FFD}[1][R]{\operatorname{FFD}(#1)}

\newcommand{\FPD}[1][R]{\operatorname{FPD}(#1)}

\newcommand{\wdt}[2][R]{\operatorname{width}_{#1}#2}
\newcommand{\dpt}[2][R]{\operatorname{depth}_{#1}#2}

\newcommand{\supP}[1]{\sup{(#1)}}
\newcommand{\infP}[1]{\inf{(#1)}}

\newcommand{\tpP}[3][R]{(\tp[#1]{#2}{#3})}

\newcommand{\DtpP}[3][R]{(\Dtp[#1]{#2}{#3})}

\newcommand{\one}{\ensuremath{\mathord \ast}}
\newcommand{\two}{\ensuremath{\mathord \dagger}}

\newcommand{\gor}{Gorenstein\xspace}
\newcommand{\mdl}{mo\-dule\xspace}              \newcommand{\mdls}[1][]{mo\-dules\xspace}
\newcommand{\cpx}{com\-plex\xspace}             
\newcommand{\mo}{mor\-phism\xspace}             
\newcommand{\ho}{ho\-mo\-\mo}                   
\newcommand{\iso}{iso\-\mo}                     
\newcommand{\qiso}{quasi-isomorphism\xspace}    
\newcommand{\dm}{du\-a\-liz\-ing \mdl}          
\newcommand{\dc}{du\-a\-liz\-ing \cpx}          
\newcommand{\sdc}{semi-\dc}                     

\renewcommand{\theequation}{\arabic{equation}}
\numberwithin{equation}{res}

\newcommand{\THom}[2]{$\operatorname{Hom}(\text{\sl #1},\text{\sl #2})$}
\newcommand{\Ttp}[2]{$\text{\sl #1}\otimes\text{\sl #2}$}

\newcommand{\fdf}[1][\f]{\operatorname{fd}#1}
\newcommand{\pdf}[1][\f]{\operatorname{pd}#1}
\newcommand{\fRS}{\mapdef{\f}{R}{S}}
\newcommand{\fRSl}{\mapdef{\f}{\Rmk}{\Snl}}

\newcommand{\crs}[1]{\pmb{#1}}
\newcommand{\s}[1]{\tilde{#1}}

\newcommand{\HLX}{(g)}
\newcommand{\HPB}{(e)}
\newcommand{\HPA}{(f)}
\newcommand{\HLA}{(f')}

\newcommand{\HBE}{(i)}
\newcommand{\HAE}{(h)}
\newcommand{\TXP}{(b)}

\newcommand{\TBL}{(c')}

\newcommand{\TAF}{(a)}
\newcommand{\TBF}{(c)}

\newcommand{\HGP}{(g)}
\newcommand{\HAEE}{(d)}
\newcommand{\HXEE}{(e)}
\newcommand{\HGE}{(e')}
\newcommand{\HXF}{(f)}
\newcommand{\HGF}{(f')}
\newcommand{\HEB}{(i)}
\newcommand{\HFB}{(h)}
\newcommand{\HFBB}{(d)}
\newcommand{\TPX}{(b)}
\newcommand{\HPBB}{(h')}

\newcommand{\TEG}{(c')}
\newcommand{\TFX}{(a')}

\newcommand{\TEA}{(c)}
\newcommand{\TFA}{(a)}

\newcommand{\pR}{f.g.\,projective/R\xspace}
\newcommand{\PR}{projective/R\xspace}
\newcommand{\IR}{injective/R\xspace}
\newcommand{\FR}{flat/R\xspace}
\newcommand{\GpR}{f.g.\,G--projective/R\xspace}
\newcommand{\GPR}{G--projective/R\xspace}
\newcommand{\GIR}{G--injective/R\xspace}
\newcommand{\GFR}{G--flat/R\xspace}
\newcommand{\PS}{projective/S\xspace}
\newcommand{\IS}{injective/S\xspace}
\newcommand{\FS}{flat/S\xspace}
\newcommand{\GPS}{G--projective/S\xspace}
\newcommand{\GIS}{G--injective/S\xspace}
\newcommand{\GFS}{G--flat/S\xspace}

\newcommand{\sym}[1]{$#1$}
\newcommand{\cit}[1]{\scriptsize #1} 
\newenvironment{stabmatrix}[1]{\begin{equation*}\label{tab:#1}}{\end{equation*}}

\newcommand{\unit}[1]{\eta_{#1}}
\newcommand{\counit}[1]{\varepsilon_{#1}}
\newcommand{\Kc}[2][R]{\operatorname{K}_{#1}[#2]}

\renewcommand{\tev}[2][RS]{\omega_{#2}^{#1}}
\renewcommand{\hev}[2][RS]{\theta_{#2}^{#1}}
\renewcommand{\A}[2][R]{\mbox{\ensuremath{_{\scriptstyle #2}\Cat{#1}{A}}}}
\renewcommand{\B}[2][R]{\mbox{\ensuremath{_{\scriptstyle #2}\mathsf{B}(#1)}}}

\begin{document}

%\sloppy

\title{Ascent Properties of Auslander Categories}

\author{Lars Winther Christensen\ \ and\ \ Henrik Holm}

\thanks{L.W.C. was partly supported by a grant from the Danish Natural
  Science Research Council.}

\address{Lars Winther Christensen, Department of Mathematics,
  University of Nebraska, Lincoln, NE 68588-0130, U.S.A.}

\urladdr{http://www.math.unl.edu/{\tiny $\sim$}lchristensen}

\email{winther@math.unl.edu}

\address{Henrik Holm, Department of Mathematical Sciences, University
  of Aarhus, Ny Munke\-gade Bldg.\ 530, DK-8000 Aarhus C, Denmark}

\urladdr{http://home.imf.au.dk/holm}

\email{holm@imf.au.dk}

%\date{\today}


%% 13D05=Homological dimension
%% 13D07 Homological functors on modules (Tor, Ext, etc.)
%% 13D25 Complexes

\begin{abstract}
  Let $R$ be a homomorphic image of a Gorenstein local ring. Recent
  work has shown that there is a bridge between Auslander categories
  and modules of finite Gorenstein homological dimensions over $R$.
  
  We use Gorenstein dimensions to prove new results about Auslander
  categories and vice versa. For example, we establish base change
  relations between the Auslander categories of the source and target
  rings in a homomorphism $\fRS$ of finite flat dimension.
\end{abstract}

\maketitle

%%% INTRODUCTION
\section*{Introduction}

\noindent
Transfer of homological properties along ring homomorphisms is already
a classical field of study, initiated in \cite{egaIV} and continued in
the more recent series
\cite{LLAHBF92,LLAHBF97,LLAHBF98,AFH-94,AFL-93}. In this paper we
investigate ascent properties of modules in the so-called Auslander
categories of a commutative noetherian ring.

For a local ring $R$ with a dualizing complex, Avramov and Foxby
\cite{LLAHBF97} introduced the Auslander categories $\A{}$ and $\B{}$,
two subcategories of the derived category of $R$. This was part of
their study of local ring homomorphisms of finite Gorenstein
dimension. One theme played in \cite{LLAHBF97} is
\begin{center}
  (I) \textit{Results for Auslander categories have implications for
    Gorenstein dimensions}
\end{center}
This is based on the realization that Auslander categories and
Gorenstein homological dimensions are close kin \cite{CFH-,EJX-96b}.
The latter were introduced much earlier by Auslander and Bridger
\cite{MAs67,MAsMBr69} and Enochs, Jenda
et.~al.~\cite{EEnOJn95b,EJT-93}.

\begin{ipg}
  In this paper we continue the theme (I). Let $\fRS$ be a local
  homomorphism of rings. Working directly with the definition of
  $\sf{A}$ we prove e.g.\ \prppartref[]{ABstabC}{c}:
  \begin{thmquote}
    \textbf{Theorem I.} \sl Assume that $\fdf$ is finite and $S$ has a
    dualizing complex. If $P$ is an $R$--module of finite projective
    dimension and $\s{A}\in\A[S]{}$ then $\DHom{P}{\s{A}}$ belongs to
    $\A[S]{}$.
  \end{thmquote}
  [Here $\operatorname{\mathbf{R}Hom}$ is the right-derived Hom
  functor.] For the next result we need the notion of Gorenstein flat
  modules, which is a generalization of flat modules introduced in
  \cite{EJT-93}. Theorem I has as a consequence \prppartref[]{HPA}{b}:
  \begin{thmquote}
    \textbf{Corollary I.} \sl Assume that $\fdf$ is finite and $S$ has
    a dualizing complex. If $P$ is a projective $R$--module and
    $\s{A}$ is a Gorenstein flat $S$--module, then $\Hom{P}{\s{A}}$ is
    Gorenstein flat over $S$.
  \end{thmquote}
  Gorenstein dimensions and Auslander categories are truly two sides of
  one coin, and the complementary theme
  \begin{center}
    (II) \textit{Results for Gorenstein dimensions have implications for
      Auslander categories}
  \end{center}
  turns out to be equally useful. For example, from the definition of
  Gorenstein flat modules we prove \lempartref[]{HXE}{a}:
  \begin{thmquote}
    \textbf{Theorem II.} \sl Assume that $\fdf$ is finite. If $\s{F}$ is a
    flat $S$--module and $A$ is a Gorenstein flat $R$--module, then
    $\tp{\s{F}}{A}$ is Gorenstein flat over $S$.
  \end{thmquote}
  From this one gets \thmpartref[]{stabD}{a}:
  \begin{thmquote}
    \textbf{Corollary II.}  \sl Assume that $\fdf$ is finite and that $R$
    and $S$ have dualizing complexes. If $\s{F}$ is an $S$--module of
    finite flat dimension and $A\in\A{}$ then $\Dtp{\s{F}}{A}$ belongs
    to $\A[S]{}$.
  \end{thmquote}
  [Here $\otimes^\mathbf{L}$ is the left-derived tensor product
  functor.]  We are not aware of any direct proof of Corollary I, i.e.\ 
  a proof that avoids Theorem I. The same remark applies to
  Corollary/Theorem II.
\end{ipg}

\begin{ipg}
  Evaluation morphisms are important tools in the study of Auslander
  categories. Indeed, Theorem I relies on the fact that the tensor
  evaluation morphism,
  \begin{equation*}
    \tag{\one}
    \mapdef{\tev{LMN}}{\tp[S]{\Hom{L}{M}}{N}}{\Hom{L}{\tp[S]{M}{N}}},
  \end{equation*}
  is invertible when $L$ is a projective $R$--module, $M$ and $N$ are
  $S$--modules, and $N$ is finitely generated. In section
  \ref{sec:applications} we give new conditions that ensure
  invertibility of evaluation morphisms; for example \obsref[]{Gtevmod}:
  \begin{thmquote}
    \textbf{Theorem III.} \sl Assume that $\fdf$ is finite. If $L$ is
    finitely generated and Gorenstein flat over $R$, $M$ is flat over
    $S$ and $N$ is injective over $S$, then $\tev{LNM}$ in (\one) is an
    isomorphism.
  \end{thmquote}
  These new isomorphisms have applications beyond the study of Auslander
  categories, e.g.\ to formulas of the Auslander--Buchsbaum type: For a
  finitely generated $R$--module $M$ of finite flat dimension, the
  classical Auslander--Buchsbaum formula
  \begin{displaymath}
    \supremum{m\in\ZZ}{\Tor{m}{k}{M}\ne0} = \dptR - \dpt{M}
  \end{displaymath}
  is a special case of \thmpartref[]{ABE}{a}:
  \begin{equation*}
    \tag{\two}
    \dpt[S]{\DtpP{N}{M}} = \dpt[S]{N} + \dpt{M} - \dptR,
  \end{equation*}
  which holds for $R$--modules $M$ of finite flat dimension and all
  $S$--modules $N$.
  
  Results like Theorem III allow us to prove that (\two) also holds
  for $R$--modules $M$ of finite Gorenstein flat dimension and
  $S$--modules $N$ of finite injective dimension.
\end{ipg}

\begin{ipg}
  As indicated by (\two), results in this paper are stated in the
  language of derived categories; we recall the basic notions in
  section \ref{sec:notation}. The prerequisites on Auslander
  categories and Gorenstein dimensions are given in section
  \ref{sec:GorAB}. Section \ref{sec:ABstability} is devoted to the
  themes (I) and (II).  In section \ref{sec:applications} we break
  to establish certain evaluation isomorphisms and then continue the
  themes of the previous section.  In section \ref{sec:ABF} we study
  formulas of the Auslander--Buchsbaum type, and in the final,
  appendix-like, section \ref{sec:cat} we catalogue the ascent
  results obtained in sections \ref{sec:ABstability} and
  \ref{sec:applications}.
\end{ipg}

%%% SECTION 0
\setcounter{section}{-1}
\section{Notation and prerequisites}
\label{sec:notation}

\noindent 
All rings in this paper are assumed to be commutative, unital and
non-zero; throughout, $R$ and $S$ denote such rings. All modules are
unitary.

\begin{bfhpg}[Complexes]
  We denote by $\C$ the category of $R$--complexes; that is, chain
  complexes of $R$--modules. We use this notation with subscripts
  $\sqsupset$, \raisebox{-0.15ex}{$\square$}, and $\sqsubset$ to
  denote the full subcategories of left- and/or right bounded
  complexes. E.g.
  \begin{displaymath}
    X \,=\, \xymatrix{ \cdots \ar[r] & X_{n+1}
      \ar[r]^-{\partial_{n+1}^X} & X_n \ar[r]^-{\partial_n^X} & X_{n-1}
      \ar[r] & \cdots, 
    }
  \end{displaymath}
  is in $\Cbr$ if and only if $X_\l =0$ for $\l \ll 0$. We use
  superscripts ${\rm I}$, ${\rm F}$, ${\rm P}$, and ${\rm fP}$ to
  indicate that the complexes in question consist of modules which
  are, respectively, injective, flat, projective, or finite (that is,
  finitely generated) projective.

  \begin{spg}
    The notation $\D$ is used for the derived category of the abelian
    category of $R$--modules; see \cite[chap.~I]{rad} or
    \cite[chap.~10]{Wei}. We use subscripts $\sqsupset$,
    \raisebox{-0.15ex}{$\square$}, and $\sqsubset$ and superscript
    ${\rm f}$ to indicate vanishing and finiteness of homology
    modules. For homological supremum and infimum of $X\in\D$ we write
    $\sup{X}$ and $\inf{X}$. Thus, $X$ is in $\Dbr$ if and only if
    $\inf{X}>-\infty$.
    
    Since $R$ is commutative, the right derived Hom, $\DHom{-}{-}$,
    and the left derived tensor product, $\Dtp{-}{-}$, are functors
    (in two variables) in $\D$.
    
    The symbol $\eq$ denotes quasi-isomorphisms in $\C$ and
    isomorphisms in $\D$.
  
    The category of complexes of $(R,S)$--bimodules is denoted
    $\C[R,S]$. We write $\D[R,S]$ for the derived category of the
    abelian category of $(R,S)$--bimodules, and we use sub- and
    superscripts on $\D[R,S]$ as we do for $\D$.
  \end{spg}
\end{bfhpg}

\begin{bfhpg}[Homological dimensions]
  We use abbreviations pd, id, and fd for projective, injective, and
  flat dimension of complexes. By $\P$, $\I$, and $\F$ we denote the
  full subcategories of $\Db$ whose objects are complexes of finite
  projective/injective/flat dimension.
\end{bfhpg}

\begin{ipg}
  The (left derived) tensor product is left-adjoint to the (right
  derived) Hom functor; this gives the adjunction isomorphism(s). This
  and other standard isomorphisms, associativity and commutativity of
  tensor products, are used freely.
  
  The, in general not invertible, evaluation morphisms shall play a
  key role in several proofs.  For later reference, we recall a
  selection of conditions under which they are invertible.
\end{ipg}

\begin{res}[Evaluation morphisms in $\sf{C}$]
  \label{res:ev}
  Let $X\in\C$, $Y \in \C[R,S]$ and $Z\in\C[S]$. Then $\Hom[S]{Y}{Z}$,
  $\Hom{X}{Y}$, and $\tp[S]{Y}{Z}$ belong to $\C[R,S]$; the canonical
  maps
  \begin{align*}
    \hev{XYZ} \colon \tp{X}{\Hom[S]{Y}{Z}} &\,\longrightarrow\,
    \Hom[S]{\Hom{X}{Y}}{Z}\quad\text{and}\\[1ex]
    \tev{XYZ} \colon \tp[S]{\Hom{X}{Y}}{Z} &\,\longrightarrow\,
    \Hom{X}{\tp[S]{Y}{Z}}
  \end{align*}
  are morphisms in $\C[R,S]$ and functorial in $X$, $Y$, and $Z$.
  
  If two of the complexes $X$, $Y$, and $Z$ are bounded, then the
  Hom~evaluation morphism $\hev{XYZ}$ is invertible under each of the
  following extra conditions:
  \begin{prt}
  \item $X\in\CfP$; or
  \item $R$ is noetherian, $X\in\Cf$, and $Z\in\CI[S]$.
  \end{prt}
  
  If two of the complexes $X$, $Y$, and $Z$ are bounded, then the
  tensor evaluation morphism $\tev{XYZ}$ is invertible under each of
  the following extra conditions:
  \begin{prt}\setcounter{prt}{2}
  \item $X\in\CfP$;
  \item $R$ is noetherian, $X\in\Cf$, and $Z\in\CF[S]$;
  \item $Z\in\CfP[S]$; or
  \item $S$ is noetherian, $Z\in\Cf[S]$, and $X\in\CP$.
  \end{prt}
\end{res}

\begin{prf}
  Conditions \prtlbl{a}--\prtlbl{d} can be traced back to
  \cite[sec.~0, 5, and 9]{hha}, \lemcite[4.4]{LLAHBF91}, and
  \thmcite[1 and 2]{DAp99c}. We have not found references for
  \prtlbl{e} and \prtlbl{f}, so we include the argument:

  \begin{spg}
    \proofoftag{e} Under the boundedness conditions the morphism
    $\tev{XYZ}$ will, in each degree, be a finite sum of evaluation
    morphisms of modules $X_h$, $Y_i$, and $Z_j$. Thus, it is
    sufficient to deal with the module case. When $Z$ is a finite
    projective module it is a direct summand in a finite free module
    $S^\beta$. By additivity of the involved functors it suffices to
    establish the isomorphism for $Z=S^\beta$, and that follows
    immediately from the commutative diagram
    \begin{equation*}
      \begin{split}
        \xymatrix{ \tp[S]{\Hom{X}{Y}}{S^\beta} \ar[d]^-{\is}
          \ar[rr]^-{\tev{XYS^\beta}} & &
          \Hom{X}{\tp[S]{Y}{S^\beta}} \ar[d]_-{\is} \\
          \Hom{X}{Y}^\beta \ar[rr]^-{\is} & & \Hom{X}{Y^\beta} }.
      \end{split}
    \end{equation*}
    
    \proofoftag{f} As above it suffices to deal with the module case and
    we may assume that $X$ is free, $X=R^{(\Lambda)}$.  Consider the
    commutative diagram
    \begin{equation*}
      \begin{split}
        \xymatrix{\tp[S]{\Hom{R^{(\Lambda)}}{Y}}{Z} \ar[d]^-{\is}
          \ar[rr]^-{\tev{R^{(\Lambda)}YZ}} & &
          \Hom{R^{(\Lambda)}}{\tp[S]{Y}{Z}} \ar[d]_-{\is} \\
          \tp[S]{Y^{\Lambda}}{Z} \ar[rr] & & (\tp[S]{Y}{Z})^{\Lambda} }
      \end{split}
    \end{equation*}
    When $S$ is noetherian and $Z$ is finite, the lower horizontal
    homomorphism is invertible by \chpcite[II, exerc.~2]{CarEil} or
    \thmcite[3.2.22]{rha}; thus $\tev{R^{(\Lambda)}YZ}$ is an
    \iso.\qedhere
  \end{spg}
\end{prf}

\enlargethispage*{1cm}
\begin{res}[Evaluation morphisms in $\sf{D}$]
  \label{res:Dev}
  Let $X\in\D$, $Y \in \D[R,S]$ and $Z\in\D[S]$. Then
  $\DHom[S]{Y}{Z}$, $\DHom{X}{Y}$, and $\Dtp[S]{Y}{Z}$ are
  representable by complexes of $(R,S)$--bimodules. The canonical
  $R$-- and $S$--linear maps,
  \begin{align*}
    \hev{XYZ} \colon \Dtp{X}{\DHom[S]{Y}{Z}} &\,\longrightarrow\,
    \DHom[S]{\DHom{X}{Y}}{Z} \quad \text{and} \\[1ex]
    \tev{XYZ} \colon \Dtp[S]{\DHom{X}{Y}}{Z} &\,\longrightarrow\,
    \DHom{X}{\Dtp[S]{Y}{Z}}
  \end{align*}
  are functorial in $X$, $Y$, and $Z$.
  
  If $R$ is noetherian, then the Hom~evaluation morphism $\hev{XYZ}$
  is invertible, provided that:
  \begin{prt}
  \item $X \in \Pf$ and $Y \in \Db[R,S]$; or
  \item $X\in\Dfb$, $Y\in\Dbl[R,S]$, and $Z\in\I[S]$.
  \end{prt}
  
  If $R$ is noetherian, then the tensor evaluation morphism
  $\tev{XYZ}$ is invertible, provided that:
  \begin{prt}\setcounter{prt}{2}
  \item $X\in\Pf$ and $Y\in\Db[R,S]$; or
  \item $X\in\Dfb$, $Y\in\Dbl[R,S]$, and $Z\in\F[S]$.
  \end{prt}
  
  If $S$ is noetherian, then the tensor evaluation morphism
  $\tev{XYZ}$ is invertible, provided that:
  \begin{prt}\setcounter{prt}{4}
  \item $Z\in\Pf[S]$ and $Y\in\Db[R,S]$; or
  \item $X\in\P$, $Y\in\Dbr[R,S]$, and $Z\in\Dfb[S]$.
  \end{prt}
\end{res}

\begin{prf}
  Conditions \prtlbl{a}--\prtlbl{d} can be traced back to
  \cite[sec.~0, 5, and 9]{hha}, \lemcite[4.4]{LLAHBF91}, and
  \thmcite[1 and 2]{DAp99c}. Parts \prtlbl{e} and \prtlbl{f} follow
  from \resref{ev}; they have similar proofs and we only write out the
  details for \prtlbl{f}:
  
  Since $S$ is noetherian $Z\in\Dfbr[S]$ has a resolution by finite
  free $S$--modules, $\CfP[S] \ni L \xre Z$. As $X \in \P$ there also
  exists a bounded projective resolution, $\CPb \ni P \xre X$. Let $Y'
  = \Tsr{i}{Y}$ be the soft truncation of $Y$ at $i=\inf{Y}$, then $Y'
  \in \Cbr[R,S]$ is isomorphic to $Y$ in $\D[R,S]$. Now the
  tensor-evaluation morphism $\tev{XYZ}$ in $\D[R,S]$ is represented
  by
  \begin{equation*}
    \tp[S]{\Hom{P}{Y'}}{L} \xra{\tev{PY'L}}
    \Hom{P}{\tp[S]{Y'}{L}},
  \end{equation*}
  and by \respartref{ev}{e} this map is an isomorphism in $\C[R,S]$.
\end{prf}

\begin{ipg}
  All results in this paper are phrased in a \emph{relative} setting,
  that is, they refer to a homomorphism $\fRS$ of rings. We refer to
  the situation $\f = 1_R$ as \emph{the absolute case.} Complexes over
  $S$ are considered as $R$--complexes with the action given by $\f$.
  
  Again, we recall for later reference the ascent properties of the
  classical homological dimensions. To distinguish $S$--modules from
  $R$--modules we mark the former with a tilde, e.g.~$\s{N}$. This
  praxis is applied whenever convenient.
\end{ipg}

\begin{res}[Ascent for modules]
  \label{res:stab}
  Let $\fRS$ be a \ho of rings. The following hold:
  \begin{itemlist}
  \item If $\s{F}$ is a flat $S$--module and $F$ a flat $R$--module,
    then $\tp{\s{F}}{F}$ is flat over $S$
  \item If $\s{P}$ is a projective $S$--module and $P$ is a projective
    $R$--module, then $\tp{\s{P}}{P}$ is projective over $S$
  \item If $\s{F}$ is a flat $S$--module and $I$ is an injective
    $R$--module, then $\Hom{\s{F}}{I}$ is injective over $S$
  \item If $F$ is a flat $R$--module and $\s{I}$ is an injective
    $S$--module, then $\Hom{F}{\s{I}}$ is injective over $S$
  \end{itemlist}
  If $S$ is noetherian, also the following hold:
  \begin{itemlist}
  \item If $\s{I}$ is an injective $S$--module and $F$ a flat
    $R$--module, then $\tp{\s{I}}{F}$ is injective over $S$
  \item If $P$ is a projective $R$--module and $\s{F}$ a flat
    $S$--module, then $\Hom{P}{\s{F}}$ is flat over $S$
  \item If $\s{I}$ is an injective $S$--module and $I$ an injective
    $R$--module, then $\Hom{\s{I}}{I}$ is flat over $S$
  \end{itemlist}
\end{res}

\begin{prf}
  All seven results are folklore and are straightforward to verify;
  see also \cite{TIs65}.  As an example, consider the penultimate one:
  The projective module $P$ is a direct summand in a free $R$--module;
  that makes $\Hom{P}{\s{F}}$ a direct summand in a product of flat
  $S$--modules and hence flat, as $S$ is noetherian. (This and the two
  neighboring results can also be proved using the evaluation
  morphisms from \resref{ev}.)
  
  Also note that the fourth is a consequence of the isomorphism,
  \begin{equation*}
    \Hom[S]{-}{\Hom{F}{\s{I}}} \is  \Hom[S]{\tp{-}{F}}{\s{I}},
  \end{equation*}
  which is an easily verified variant of standard adjointness.
\end{prf}

\begin{bfhpg}[Ascent for complexes]
  \label{res:stabcpx}
  Let $\fRS$ be a \ho of rings. The results in \resref{stab} imply
  similar ascent results for complexes of finite homological
  dimension. In short and suggestive notation we write them as:
  \begin{itemlist}
  \item $\Dtp{\F[S]}{\F} \subseteq \F[S]$
  \item $\Dtp{\P[S]}{\P} \subseteq \P[S]$
  \item $\DHom{\F[S]}{\I} \subseteq \I[S]$
  \item $\DHom{\P}{\I[S]} \subseteq \I[S]$
  \end{itemlist}
  If $S$ is noetherian, we also have:
  \begin{itemlist}
  \item $\Dtp{\I[S]}{\F} \subseteq \I[S]$
  \item $\DHom{\P}{\F[S]} \subseteq \F[S]$
  \item $\DHom{\I[S]}{\I} \subseteq \F[S]$
  \end{itemlist}
\end{bfhpg}

\begin{bfhpg}[Local rings and homomorphisms]
  We say that $\Rmk$ is local, if $R$ is noetherian and local with
  maximal ideal $\m$ and residue field $k$. A \ho of rings $\fRSl$ is
  said to be local if $\f(\m) \subseteq \n$.
\end{bfhpg}

\begin{bfhpg}[Homomological dimensions of homomorphisms]
  \label{hdf}
  Let $\fRS$ be a homomorphism of rings. The flat dimension of $\f$ is
  by definition the flat dimension of $S$ considered as a module over
  $R$ with the action given by $\f$. That is, $\fdf = \fd{S}$. The
  projective and injective dimensions of $\f$ are defined similarly.
  
  Note that in the case where $\f$ is a local homomorphism, our
  definition of $\pd[]{\f}$ differs from the one in
  \cite[def.~4.2]{SInSSW04}. However $\pd[]{\f}$ in our sense, and
  $\pd[]{\f}$ in the sense of \cite[def.~4.2]{SInSSW04} are
  simultaneously finite.

  \begin{spg}
    Transfer of homological properties along homomorphisms is already
    a classical field of study. A basic observation is
    \prpcite[(4.6)(b)]{LLAHBF97}: If $\fdf$ is finite, then
    \begin{displaymath}
      \big(\P[S] \subseteq\!\big)\, \F[S] \subseteq \F 
      \quad \text{ and } \quad  
      \I[S] \subseteq \I.
    \end{displaymath}
    In \lemref[]{stabcpxnoeth} below we use this to establish a useful
    variant of \resref{stabcpx}.
  \end{spg}  

  \begin{spg}
    The literature emphasizes the study of homomorphisms of finite
    flat dimension; largely, we follow this tradition, as it is well-
    founded: Let $\fRS$ be a homomorphism of noetherian rings.
    \begin{itemlist}
    \item If the projective dimension of $\f$ is finite then so is
      $\fdf$, and the converse holds if $R$ has finite Krull
      dimension; see~\prpcite[6]{CUJ70} and
      \thmcite[II.3.2.6]{LGrMRn71}.
    \item If $\f$ is local, then the injective dimension of $\f$ is
      finite if and only if $\fdf$ is finite and $R$ is Gorenstein.
      This was only established recently, in \thmcite[13.2]{AIM-},
      though the surjective case goes back to
      \cite[(II.5.5)]{CPsLSz73}.
    \item If $R$ has finite Krull dimension, then $\f$ has finite
      injective dimension if and only if $\fdf$ is finite and $R$ is
      Gorenstein at any contraction $\p=\q\cap R$ of a prime ideal in
      $S$. This follows from the local case above.
    \end{itemlist}
  \end{spg}
\end{bfhpg}

\begin{ipg}
  The next lemma resembles the last part of \resref{stabcpx}; the
  difference is that the noetherian assumption has been moved from $S$
  to $R$, while the complexes all have $S$--structures.
\end{ipg}

\begin{lem}
  \label{lem:stabcpxnoeth}
  Let $\fRS$ be a \ho of rings with $\fdf$ finite. If $R$ is
  noetherian, then the following hold:
  \begin{itemlist}
  \item $\Dtp[S]{\I[S]}{\F[S]} \subseteq \I$
  \item $\DHom[S]{\P[S]}{\F[S]} \subseteq \F$
  \item $\DHom[S]{\I[S]}{\I[S]} \subseteq \F$
  \end{itemlist}
\end{lem}

\begin{prf}
  For the last assertion let $\s{I},\s{J}\in\I[S]$. We must show that
  for any finite $R$--module $M$, the complex
  $\Dtp{M}{\Hom[S]{\s{J}}{\s{I}}}$ is in $\Db$. By \respartref{Dev}{b}
  we have
  \begin{displaymath}
    \Dtp{M}{\DHom[S]{\s{J}}{\s{I}}} \eq \DHom[S]{\DHom{M}{\s{J}}}{\s{I}},
  \end{displaymath}
  and the latter complex is bounded as $\s{I} \in \I[S]$ and $\s{J}
  \in \I[S] \subseteq \I$ by \pgref{hdf}.
  
  The other assertions have similar proofs.
\end{prf}

%%% SECTION 1

\section{Gorenstein dimensions and Auslander categories}
\label{sec:GorAB}

\noindent
This paper pivots on the interplay between (semi-)dualizing complexes,
their Auslander categories, and Gorenstein homological dimensions.

Semi-dualizing complexes and Auslander categories came up in studies
of ring homomorphisms \cite{LLAHBF97} and are used to detect the
Gorenstein \cite{LWC01a} and Cohen--Macaulay \cite{HHlPJr2}
properties of rings. This section recaps the relevant definitions and
results.

\begin{bfhpg}[Gorenstein dimensions]
  \label{GorDim}
  Gorenstein projective, injective and flat modules are defined in
  terms of so-called complete resolutions:
  \begin{itemlist}
  \item An $R$--module $A$ is \emph{Gorenstein projective} if there
    exists an exact complex $\crs{P}$ of projective modules, such that
    $A \is \Coker{(P_1 \to P_0)}$ and $\H{\Hom{\crs{P}}{Q}}=0$ for all
    projective $R$--modules $Q$. Such a complex $\crs{P}$ is called a
    \emph{complete projective resolution} (of $A$).
    
  \item An $R$--module $B$ is \emph{Gorenstein injective} if there
    exists an exact complex $\crs{I}$ of injective modules, such that
    $B \is \Ker{(I_0 \to I_{-1})}$ and $\H{\Hom{J}{\crs{I}}}=0$ for
    all injective $R$--modules $J$. Such a complex $\crs{I}$ is called
    a \emph{complete injective resolution} (of $B$).
    
  \item An $R$--module $A$ is \emph{Gorenstein flat} if there exists
    an exact complex $\crs{F}$ of flat modules, such that $A \is
    \Coker{(F_1 \to F_0)}$ and $\H{\tp{J}{\crs{F}}}=0$ for all
    injective $R$--modules $J$. Such a complex $\crs{F}$ is called a
    \emph{complete flat resolution} (of $A$).
  \end{itemlist}
  These definitions from \cite{EEnOJn95b, EJT-93} generalize and
  dualize the notion of G--dimension 0 modules from
  \cite{MAs67,MAsMBr69}; see \thmcite[(4.2.5) and (5.1.11)]{LWC}.

  \begin{spg}
    By taking resolutions, one defines the \emph{Gorenstein
      projective} and \emph{Gorenstein flat dimension} of
    right-bounded complexes in $\D$ and \emph{Gorenstein injective
      dimension} of left-bounded complexes. For details see
    \cite[def.~(4.4.2), (5.2.2), and (6.2.2)]{LWC}. All projective
    modules are Gorenstein projective, so the Gorenstein projective
    dimension of an $R$--complex $X\in\Dbr$ is a finer invariant
    than the usual projective dimension; that is, $\Gpd{X} \le
    \pd{X}$. Similarly, injective and flat modules are Gorenstein
    injective and Gorenstein flat, so we have inclusions
    \begin{equation*}
      \P \subseteq \GP, \ \ \I \subseteq \GI \ \ \text{ and } \ \
      \F \subseteq \GF.  
    \end{equation*}
    Here $\GP$ denotes the full subcategory of bounded complexes of
    finite Gorenstein projective dimension; $\GI$ and $\GF$ are defined
    similarly. See \cite{HHl04a, LWC,CFH-} for details on Gorenstein
    dimensions.
  \end{spg}
\end{bfhpg}

\begin{bfhpg}[Auslander Categories]
  \label{res:ac}
  Assume that $R$ is noetherian. A \emph{semi-dualizing complex} for
  $R$ is a complex $C \in \Dfb$ such that the homothety morphism,
  \begin{equation*}
    R \longrightarrow \DHom{C}{C},
  \end{equation*}
  is invertible in $\D$; cf.~\cite{LWC01a}. Note that $R$ is a
  semi-dualizing complex for itself.
  
  If, in addition, $C \in \I$, then $C$ is a \emph{dualizing complex}
  for $R$, cf.~\cite[V.\S2]{rad}.

  \begin{spg}
    Let $C$ be a semi-dualizing complex for $R$, and consider the
    adjoint pair of functors
    \begin{displaymath}
      \tag{\one}
      \xymatrix{\D \ar@<0.7ex>[rrr]^-{C\otimes_R^{\bf L}-} & {} & {} & 
        \D. \ar@<0.7ex>[lll]^-{\DHom{C}{-}} 
      } 
    \end{displaymath}
    The Auslander categories with respect to $C$, denoted $\A{C}$ and
    $\B{C}$, are the full subcategories of $\Db$ whose objects are
    specified as follows:
    \begin{align*}
      \A{C} &\,=\, \left\{ X \in \Db \: \left|
          \begin{array}{l}
            \mbox{ $\eta_X \colon X \xre \DHom{C}{\Dtp{C}{X}}$ is an iso-} \\
            \mbox{ morphism in $\D$, and $\Dtp{C}{X} \in \Db$} 
          \end{array}
        \right.
      \right\}, \\ \\
      \B{C} & \,=\, \left\{ Y \in \Db \: \left|
          \begin{array}{l}
            \mbox{ $\varepsilon_Y \colon \Dtp{C}{\DHom{C}{Y}}\xre Y$ is an
              isomor-} \\ 
            \mbox{ phism in $\D$, and $\DHom{C}{Y} \in \Db$}
          \end{array}
        \right.  \right\},
    \end{align*}
    where $\unit{}$ and $\counit{}$ denote the unit and counit of the
    pair $(\Dtp{C}{-}, \DHom{C}{-})$. These categories were introduced
    in \cite{LLAHBF97,LWC01a}.
    
    The Auslander categories are triangulated subcategories of $\D$,
    and the adjoint pair in (\one) restricts to an equivalence between
    them,
    \begin{displaymath}
      \xymatrix{\A{C}
        \ar@<0.7ex>[rrr]^-{C\otimes_R^{\mathrm{\mathbf{L}}}-} 
        & {} & {} &
        \B{C}. \ar@<0.7ex>[lll]^-{\DHom{C}{-}}
      }
    \end{displaymath}
    By \cite[prop.~(4.4)]{LWC01a} there are inclusions $\F \subseteq
    \A{C}$ and $\I \subseteq \B{C}$.
  \end{spg}

  \begin{spg}
    The relation between Auslander categories and Gorenstein dimensions
    is established in \cite{EJX-96b,SYs95b,CFH-}: If $D$ is a dualizing
    complex for $R$, then
    \begin{equation}
      \label{eq:main}
      \A{D} \,=\, \GP \,=\, \GF \quad \text{ and } \quad \B{D} \,=\,
      \GI. 
    \end{equation}
  \end{spg}
\end{bfhpg}

%%% SECTION 2

\section{Ascent properties}
\label{sec:ABstability}

\noindent
The first result below should be compared to \resref{stabcpx}.

\begin{prp}
  \label{prp:ABstabC}
  Let $\mapdef{\f}{R}{S}$ be a \ho of rings. If $S$ is noetherian and
  $\s{C}$ is a \sdc for $S$, then the following hold:
  \begin{prt}
  \item If $\s{A}\in\A[S]{\s{C}}$ and $F\in\F$, then $\Dtp{\s{A}}{F}
    \in \A[S]{\s{C}}$
  \item If $\s{B}\in\B[S]{\s{C}}$ and $F\in\F$, then $\Dtp{\s{B}}{F}
    \in \B[S]{\s{C}}$
  \item If $P\in\P$ and $\s{A}\in\A[S]{\s{C}}$, then $\DHom{P}{\s{A}}
    \in \A[S]{\s{C}}$
  \item If $P\in\P$ and $\s{B}\in\B[S]{\s{C}}$, then $\DHom{P}{\s{B}}
    \in \B[S]{\s{C}}$
  \item If $\s{A}\in\A[S]{\s{C}}$ and $I\in\I$, then $\DHom{\s{A}}{I}
    \in \B[S]{\s{C}}$
  \item If $\s{B}\in\B[S]{\s{C}}$ and $I\in\I$, then $\DHom{\s{B}}{I}
    \in \A[S]{\s{C}}$
  \end{prt}
  \begin{spg}
    For the important special case of a dualizing complex, absolute
    versions of parts \prtlbl{a}, \prtlbl{b}, \prtlbl{e}, and \prtlbl{f}
    appear in \cite[(6.4.13)]{LWC} and certain relative versions in
    \cite{LKhSYs03}.
  \end{spg}
\end{prp}

\begin{prf}
  \proofoftag{c} First note that $\DHom{P}{\s{A}}$ belongs to
  $\Db[S]$, as $P$ has finite projective dimension over $R$ and
  $\s{A}\in\Db[S]$. To see that also $\Dtp[S]{\s{C}}{\DHom{P}{\s{A}}}$
  is homologically bounded, we employ \respartref{Dev}{f} to get an
  isomorphism,
  \begin{equation*}
    \Dtp[S]{\s{C}}{\DHom{P}{\s{A}}} \xra[\eq]{\quad\omega\quad}
    \DHom{P}{\Dtp[S]{\s{C}}{\s{A}}}.
  \end{equation*}
  The latter complex is homologically bounded as
  $\Dtp[S]{\s{C}}{\s{A}}$ is so. Finally, the commutative diagram
  \begin{equation*}
    \xymatrix{\DHom{P}{\s{A}} \ar[d]^-{\eq}_-{\DHom{P}{\unit{\s{A}}}}
      \ar[rr]^-{\unit{\DHom{P}{\s{A}}}} & &
      \DHom[S]{\s{C}}{\Dtp[S]{\s{C}}{\DHom{P}{\s{A}}}}
      \ar[d]_-{\eq}^-{\DHom[S]{\s{C}}{\omega}} \\
      \DHom{P}{\DHom[S]{\s{C}}{\Dtp[S]{\s{C}}{\s{A}}}}
      \ar[rr]^-{\eq}_-{\text{swap}} & & 
      \DHom[S]{\s{C}}{\DHom{P}{\Dtp[S]{\s{C}}{\s{A}}}} }
  \end{equation*}
  shows that the unit $\unit{\DHom{P}{\s{A}}}$ is invertible in $\D$.
  
  The proof of \prtlbl{d} is similar and also uses
  \respartref{Dev}{f}; the proofs of \prtlbl{a}, \prtlbl{b},
  \prtlbl{e}, and \prtlbl{f} are also similar and rely on standard
  conditions for invertibility of tensor and Hom~evaluation morphisms,
  cf.~\resref{Dev}.
\end{prf}

\begin{ipg}
  The next example shows why \resref{stabcpx} and \prpref{ABstabC}
  give no results about $\Dtp{\I}{\I}$ and $\DHom{\I}{\P}$.
\end{ipg}

\begin{exa}
  \label{exa:counter}
  Let $k$ be a field and consider the zero-dimensional local ring
  \begin{displaymath}
    R \,=\, \pows{X,Y}/(X^2\!,XY,Y^2). 
  \end{displaymath}
  We write the residue classes $x=[X]$ and $y=[Y]$; since $(x) \cap
  (y) = (0)$, the ring $R$ is not \gor, but it has an injective \dm
  $D=\E{k}$.  Furthermore, the maximal ideal $(x,y)$ is isomorphic to
  $k^2$, so we have an exact sequence $0 \to k^2 \to R \to k \to 0$.
  Applying $\Hom{-}{D}$ we get:
  \begin{align*}
    \tag{\text{$*$}} 0 \to k \to D \to k^2 \to 0.
  \end{align*}
  
  \begin{spg}
    \prtlbl{a} Since $R$ is not \gor, $D\notin\A{D}$ by
    \thmcite[(3.3.5)]{LWC}, and hence $\Dtp{D}{D}$ is not in $\B{D}$;
    cf.~\thmcite[(3.2)]{LLAHBF97}. Of course, this argument is valid
    over any non-\gor local ring with a \dc; over the ring in question
    it is even true that $\Dtp{D}{D}$ is not homological bounded, in
    particular $\Dtp{D}{D}$ is not in $\A{D}$ either. To see this, we
    assume that $\Tor{m}{D}{D}=0$ for $m \ge$ some $m_0$ and use $(*)$
    to derive a contradiction.  Applying $\Dtp{D}{-}$ to $(*)$ we get
    a long exact sequence of Tor-modules, which shows that
    $\Tor{m}{D}{k} \is \Tor{m+1}{D}{k^2}$ for $m \ge m_0$. Thus, the
    Betti numbers of $D$ satisfy the relation $\bett{m_0}{D} =
    2^t\bett{m_0+t}{D}$ for $t\ge 0$. But this implies that they must
    vanish from $m_0$, in particular $D$ must have finite projective
    dimension, which is tantamount to $R$ being \gor.
  \end{spg}

  \begin{spg}
    \prtlbl{b} Since $R$ is not \gor, $\DHom{D}{R}$ does not belong to
    $\A{D}$ cf.~\thmcite[(3.2)]{LLAHBF97} and \thmcite[(3.3.5)]{LWC}.
    Again, this argument is valid for any non-\gor local ring with a
    \dc. For this specific ring the complex $\DHom{D}{R}$ is not even
    bounded, and in particular not in $\B{D}$. One can show this by
    applying $\DHom{-}{R}$ to $(*)$ and then argue as in (a).
  \end{spg}

  \begin{spg}
    These arguments actually show that any finite module in $\A{D}$
    has finite flat dimension, and any finite module in $\B{D}$ has
    finite injective dimension. The first part, at least, is
    well-known as $R$ is Golod and not a hypersurface,
    cf.~\cite[exa.~3.5(2)]{LLAAMr02}. However, there is a direct
    argument that applies to arbitrary modules; it is a non-finite
    version of \prpcite[2.4]{YYs03} and equivalent to
    \prpcite[6.1(2)]{SInHKr}:
  \end{spg}
\end{exa}

\begin{prp}
  Let $\Rmk$ be local with $\m^2 =0$. If $R$ is not Gorenstein, then
  any Gorenstein flat $R$--module is free, and any Gorenstein
  injective $R$--module is injective.
\end{prp}

\begin{prf}
  Let $A$ be Gorenstein flat.  Since the maximal ideal of $R$ is
  nilpotent \prpcite[3 and 15]{SEl56} gives a projective cover of $A$,
  that is, an exact sequence
  \begin{equation*}
    \tag{\one}
    0 \to K \to P \to A \to 0,
  \end{equation*}
  where $P$ is projective and $K \subseteq \m P$. The module $K$ is a
  $k$--vector space, because $\m^2 = 0$, and Gorenstein flat by
  exactness of (\one). If $K\ne 0$ this implies that $k$ is Gorenstein
  flat, which contradicts the assumption that $R$ is not Gorenstein.
  Hence, $K=0$ and $A$ is isomorphic to $P$, which is free as $R$
  is local, cf.~\thmcite[2]{IKp58}.
  
  If $B$ is Gorenstein injective, $\Hom{B}{\E{k}}$ is Gorenstein flat
  by \prpcite[5.1]{CFH-} and hence flat. Since $\E{k}$ is faithful,
  this implies that $B$ is injective.
\end{prf}

\begin{ipg}
  Parts \prtlbl{a} and \prtlbl{c} in the next proposition are relative
  versions of \corcite[5.2 and prop.~5.1]{CFH-}. The proofs presented
  here require existence of a dualizing complex $\s{D}$ for $S$ so
  that results about Auslander categories from \prpref{ABstabC} may be
  applied to $\s{C}=\s{D}$. We do not know of a proof that does not
  use dualizing complexes.
\end{ipg}

\begin{prp}
  \label{prp:HPA}
  Let $\fRS$ be a \ho of rings with $\fdf$ finite. If $S$ is
  noetherian and admits a \dc, then the following hold:
  \begin{prt}
  \item If $\s{B}$ is \gor injective over $S$\ and $F$ is flat over
    $R$, then $\tp{\s{B}}{F}$ is \gor injective over $S$
  \item If $P$ is projective over $R$ and $\s{A}$ is \gor flat over
    $S$, then $\Hom{P}{\s{A}}$ is \gor flat over $S$
  \item If $\s{B}$ is \gor injective over $S$\ and $I$ is injective
    over $R$, then $\Hom{B}{I}$ is \gor flat over $S$
  \end{prt}
\end{prp}

\begin{prf}
  The three assertions have similar proofs; we only write out part
  \prtlbl{b}: Let $\s{D}$ be a dualizing complex for $S$. The module
  $\Hom{P}{\s{A}}$ represents $\DHom{P}{\s{A}}$, so
  $\Gfd[S]{\Hom{P}{\s{A}}}$ is finite by \prppartref{ABstabC}{c} and
  \eqref{main}.  Actually,
  \begin{equation*}
    \Gfd[S]{\Hom{P}{\s{A}'}} \le d := \FFD[S]<\infty
  \end{equation*}
  for any Gorenstein flat $S$--module $\s{A}'$ by \thmcite[3.5]{CFH-}
  and \corcite[V.7.2]{rad}. Here $\FFD[S]$ is the finitistic flat
  dimension of $S$, which is defined as
  \begin{equation*}
    \FFD[S] = \sup\{\fd[S]{\s{M}} \mid \text{$\s{M}$ is an $S$--module
      of finite flat dimension}\}.
  \end{equation*}
  Consider a piece of a complete flat resolution of $\s{A}$:
  \begin{equation*}
    0 \to \s{A} \to \s{F}_0 \to \cdots \to \s{F}_{d-1} \to \s{A}' \to 0;
  \end{equation*}
  also $\s{A}'$ is Gorenstein flat. Applying the exact functor
  $\Hom{P}{-}$ gives an exact sequence of $S$--modules,
  \begin{equation*}
    0 \to \Hom{P}{\s{A}} \to \Hom{P}{\s{F}_0} \to \cdots \to
    \Hom{P}{\s{F}_{d-1}} \to \Hom{P}{\s{A}'} \to 0.
  \end{equation*}
  By \resref{stab} the modules $\Hom{P}{\s{F}_\l}$ are flat over $S$,
  whence $\Hom{P}{\s{A}}$ is Gorenstein flat over $S$ by
  \thmcite[3.14]{HHl04a}.
\end{prf}

\begin{rmk}
  \prpref{HPA} demonstrates how ascent properties of Auslander
  categories yield ascent results for Gorenstein dimensions, and we do
  not know any other way to prove these results.  However, information
  also flows in the opposite direction. The following \lemref{HXE}
  about Gorenstein dimensions is the corner stone in the proof of the
  subsequent \thmref{stabD} concerning Auslander categories.  Again,
  we are not aware of any proof of \thmref[]{stabD} that does not use
  the connection to Gorenstein dimensions.
\end{rmk}

\begin{lem} 
  \label{lem:HXE} 
  Let $\fRS$ be a \ho of rings with $\fdf$ finite, then the following
  hold:
  \begin{prt}
  \item If $\s{F}$ is flat over $S$ and $A$ is \gor flat over $R$,
    then $\tp{\s{F}}{A}$ is \gor flat over $S$.
  \item If $\s{P}$ is projective over $S$ and $B$ is \gor injective
    over $R$, then $\Hom{\s{P}}{B}$ is \gor injective over $S$,
    provided that $\P[S] \subseteq \P$.
  \item If $A$ is \gor flat over $R$ and $\s{I}$ is injective over
    $S$, then $\Hom{A}{\s{I}}$ is \gor injective over $S$.
  \end{prt}
\end{lem}

\begin{prf}
  \proofoftag{a} Let $\crs{Q}$ be a complete flat resolution of $A$.
  The module $\s{F}$ has finite flat dimension over $R$, so the
  complex $\tp{\s{F}}{\crs{Q}}$ of flat $S$--modules,
  cf.~\resref{stab}, is exact by \lemcite[2.3]{CFH-}.  For any
  injective $S$--module $\s{J}$ we have exactness of
  $\tp[S]{\s{J}}{\tpP{\s{F}}{\crs{Q}}} \is
  \tp[S]{\s{F}}{\tpP{\s{J}}{\crs{Q}}}$, as $\s{J}\in\I[S] \subseteq
  \I$, so $\tp{\s{F}}{\crs{Q}}$ is a complete flat resolution over
  $S$.
  
  \proofoftag{b} Let $\crs{I}$ a complete injective resolution of $B$,
  then $\Hom{P}{\crs{I}}$ is a complex of injective $S$--modules,
  cf.~\resref{stab}, and exact as $\s{P}\in\P[S] \subseteq \P$. For
  any injective $S$--module $\s{J}$, the complex
  $\Hom[S]{\s{J}}{\Hom{\s{P}}{\crs{I}}} \is
  \Hom[S]{\s{P}}{\Hom{\s{J}}{\crs{I}}}$ is exact, as $\s{J}\in \I[S]
  \subseteq \I$.
  
  \proofoftag{c} Let $\crs{F}$ be a complete flat resolution of $A$.
  By \resref{stab} the complex $\Hom{\crs{F}}{\s{I}}$ consists of
  injective $S$--modules. To see that it is exact, write
  $\Hom{\crs{F}}{\s{I}} \is \Hom[S]{\tp{\crs{F}}{S}}{\s{I}}$.
  Exactness then follows as $S \in \F$; cf.~\lemcite[2.3]{CFH-}.  For
  any injective $S$--module $\s{J}$, we have
  $\Hom[S]{\s{J}}{\Hom{\crs{F}}{\s{I}}} \is
  \Hom[S]{\tp{\s{J}}{\crs{F}}}{\s{I}}$, which is exact as $\s{J}\in
  \I[S] \subseteq \I$.
\end{prf}

\begin{bfhpg}[Ascent results for complexes]
  \label{res:SRFC}
  Just like the classical ascent results for modules, \resref{stab},
  the lemma above gives rise to ascent results for complexes, similar
  to \resref{stabcpx}. For example, if $\fRS$ is of finite flat
  dimension, then
  \begin{equation*}
    \tag{a}
    \Dtp{\F[S]}{\GF} \subseteq \GF[S].
  \end{equation*}
  Indeed, let $\s{F}$ and $A$ be bounded complexes of, respectively,
  flat $S$--modules and \gor flat $R$--modules. Since the modules in
  $\s{F}$ have finite flat dimension over $R$, the bounded complex
  $\tp{\s{F}}{A}$ represents $\Dtp{\s{F}}{A}$ by \corcite[2.16]{CFH-}.
  \lempartref{HXE}{a} --- and the fact that direct sums of Gorenstein
  flat modules are Gorenstein flat \prpcite[3.2]{HHl04a} --- now shows
  that $\tp{\s{F}}{A}$ is a complex of \gor flat $S$--modules; thus it
  belongs to $\GF[S]$. We even get a bound on the dimension:
  \begin{equation*}
    \tag{b}
    \Gfd[S]{\DtpP{\s{F}}{A}} \le \fd[S]{\s{F}} + \Gfd{A} \quad\text{for
      $\s{F}\in\F[S]$ and $A\in\GF$}.
  \end{equation*}
  A series of ascent results for modules --- including \prpref[]{HPA}
  and \lemref[]{HXE} --- are summed up in tables \pgref{table1} and
  \pgref{table2}. We do not write out the corresponding results for
  complexes, let alone the bounds on homological dimensions. Rather,
  we leave it to the reader to derive further results like (a) and (b)
  above from the module versions in \pgref{table1} and \pgref{table2}.
  Certain results of this kind have been established under more
  restrictive conditions in \seccite[6.4]{LWC} and
  \seccite[3]{LKhSYs03}.
\end{bfhpg}

\begin{thm}
  \label{thm:stabD}
  Let $\fRS$ be a \ho of noetherian rings with $\fdf$ finite. If $R$
  and $S$ have dualizing complexes $D$ and $\s{D}$, respectively, then
  the following hold:
  \begin{prt}
  \item If $\s{F}\in\F[S]$ and $A\in\A{D}$, then $\Dtp{\s{F}}{A} \in
    \A[S]{\s{D}}$
  \item If $\s{P}\in\P[S]$ and $B\in\B{D}$, then $\DHom{\s{P}}{B} \in
    \B[S]{\s{D}}$
  \item If $A\in\A{D}$ and $\s{I}\in\I[S]$, then $\DHom{A}{\s{I}} \in
    \B[S]{\s{D}}$
  \end{prt}  
\end{thm}

\begin{prf}
  Recall that in the presence of dualizing complexes, we have
  $\A{D}=\GF$ and $\A[S]{\s{D}}=\GF[S]$. Part \prtlbl{a} is now a
  reformulation of \respartref{SRFC}{a} above. Also \prtlbl{b} and
  \prtlbl{c} are straightforward consequences of \lemref{HXE}. For
  part \prtlbl{b} note that $\P[S] \subseteq \F$ since $\fdf$ is
  finite, and furthermore $\F = \P$ as $R$ has a dualizing complex;
  see e.g.\ \cite[proof of cor.~3.4]{HBF77b}.
\end{prf}

\begin{ipg}
  We should like to stress that no assumptions are made in
  \thmref{stabD} about an explicit connection between the dualizing
  complexes $D$ and $\s{D}$. If, on the other hand, $\s{D}$ is the
  base change of $D$, then it is elementary to verify \thmref[]{stabD}
  from the definitions of Auslander categories; see e.g.\ 
  \prpcite[(5.8)(a)]{LWC01a}.
\end{ipg}

%%% SECTION 3

\section{Evaluation morphisms}
\label{sec:applications}

\noindent
The main results of this section --- \thmref[Theorems~]{Ghev} and
\thmref[]{Gtev} --- give new sufficient conditions for invertibility
of evaluation morphisms. To get a feeling for the nature of these
results, compare \thmpartref{Ghev}{a} below to \respartref{Dev}{a}:
The condition on the left-hand complex, $X$, has been relaxed from
finite projective dimension to finite Gorenstein projective dimension,
and in return conditions of finite homological dimension have been
imposed on the other two complexes.

However special the conditions in \thmref[]{Ghev} and \thmref[]{Gtev}
may seem, the theorems have interesting applications; these are
explored in \prpref[]{stabDS}--\corref[]{HXF} below and further in
section \ref{sec:ABF}. The proofs of \thmref[]{Ghev} and
\thmref[]{Gtev} are deferred to the end of the section.

\begin{thm}[Hom evaluation]
  \label{thm:Ghev}
  Let $\fRS$ be a homomorphism of rings with $\fdf$ finite.  For
  complexes $X\in\D$ and $Y,Z\in\D[S]$ the Hom~evaluation morphism
  \begin{equation*}
    \dmapdef{\hev{XYZ}}{\Dtp{X}{\DHom[S]{Y}{Z}}}{\DHom[S]{\DHom{X}{Y}}{Z}}
  \end{equation*}
  is an isomorphism in $\D[S]$, provided that
  \begin{prt}
  \item $X\in\GPf$, $Y\in\P[S]$, $Z\in\F[S]$, and $R$ is noetherian.
  \end{prt}
  
  For complexes $U,V\in\D[S]$ and $W\in\D$ the \mo
  \begin{equation*}
    \dmapdef{\hev[SR]{UVW}}{\Dtp[S]{U}{\DHom{V}{W}}}{\DHom{\DHom[S]{U}{V}}{W}}
  \end{equation*}
  is an isomorphism in $\D[S]$ provided that
  \begin{prt} \setcounter{prt}{1}
  \item $U\in\If[S]$, $V\in\I[S]$, $W\in\GI$, $S$ is noetherian, and
    $\F[S] \subseteq \P$\footnote{\label{fnt}In
      \thmref[Theorems~]{Ghev} and \thmref[]{Gtev} we encounter two
      requirements:
      \begin{displaymath}
        \P[S] \subseteq \P \quad \textnormal{ and } \quad \F[S] \subseteq
        \P; 
      \end{displaymath}
      in fact we already met the former in \lemref{HXE}. In the absolute
      case, the first one is void while the second says that flat modules
      have finite projective dimension.
      
      It is clear that the first is weaker than the second and tantamount
      to $\pdf$ being finite.  By \prpcite[6]{CUJ70}, the second
      requirement is satisfied when:
      \begin{rqm}
      \item $\fdf$ is finite and $\FPD$ is finite; or
      \item $\pdf$ is finite and $\FPD[S]$ is finite.
      \end{rqm}
      Here $\FPD$ is the finitistic projective dimension of $R$, which is
      defined as
      \begin{equation*}
        \FPD = \sup \{\pd{M} \mid \text{$M$ is an $R$--module of finite
          projective dimension}\}.
      \end{equation*}
      Recall that over a noetherian ring, the finitistic projective
      dimension, $\operatorname{FPD}$, is equal to the Krull dimension;
      see \corcite[5.5]{HBs62} and \thmcite[II.3.2.6]{LGrMRn71}.
      
      The conditions (1) and (2) are actually independent: Let $Q$
      denote Nagata's noetherian regular ring of infinite Krull
      dimension \cite[example~1, p.~203]{Nag} and consider the natural
      inclusion $k \into Q$, where $k$ is the field over which $Q$ is
      built.  Clearly, $k \into Q$ satisfies (1) but not (2).  Now let
      $\m$ be any maximal ideal of $Q$ and consider the projection $Q
      \onto Q/\m$.  Since $Q$ is regular, every finite $Q$--module has
      finite projective dimension, cf.~\corcite[3]{SGt82}, ergo $Q
      \onto Q/\m$ satisfies (2) but not (1).}.
  \end{prt}
\end{thm}

\begin{thm}[Tensor evaluation]
  \label{thm:Gtev}
  Let $\fRS$ be a homomorphism of rings with $\fdf$ finite.  For
  complexes $U,V \in\D[S]$ and $W\in\D$ the tensor evaluation morphism
  \begin{equation*}
    \dmapdef{\tev[SR]{UVW}}{\Dtp{\DHom[S]{U}{V}}{W}}{\DHom[S]{U}{\Dtp{V}{W}}}
  \end{equation*}
  is an isomorphism in $\D[S]$ if either:
  \begin{prt}
  \item $U\in\If[S]$, $V\in\I[S]$, $W\in\GF$, and $S$ is noetherian;
    or
  \item $U\in\I[S]$, $V\in\I[S]$, $W\in\GFf$, and $R$ is noetherian.
  \end{prt}
  
  For complexes $X\in\D$ and $Y,Z\in\D[S]$ the morphism
  \begin{equation*}
    \dmapdef{\tev{XYZ}}{\Dtp[S]{\DHom{X}{Y}}{Z}}{\DHom{X}{\Dtp[S]{Y}{Z}}}
  \end{equation*}
  is an isomorphism in $\D[S]$ if:
  \begin{prt} \setcounter{prt}{2}
  \item $X\in\GPf$, $Y\in\F[S]$, $Z\in\I[S]$, and $R$ is noetherian;
  \item $X\in\GP$, $Y\in\F[S]$, $Z\in\If[S]$, $S$ is noetherian, and
    $\F[S]\subseteq \P^{\ref{fnt}}$; or
  \item[\prevprt{'}] $X\in\GP$, $Y\in\P[S]$, $Z\in\If[S]$, $S$ is
    noetherian, and $\P[S]\subseteq \P^{\ref{fnt}}$.
  \end{prt}
\end{thm}

\begin{obs}
  \label{obs:Gtevmod}
  Let $R$ be noetherian and $G$ be a finite Gorenstein projective
  $R$--module; let $\s{F}$ be a flat $S$--module and $\s{I}$ an
  injective $S$--module. Then $\Dtp[S]{\DHom{G}{\s{F}}}{\s{I}}$ is
  represented by $\tp[S]{\Hom{G}{\s{F}}}{\s{I}}$ and
  $\DHom{G}{\Dtp[S]{\s{F}}{\s{I}}}$ by
  $\Hom{G}{\tp[S]{\s{F}}{\s{I}}}$; thus \thmpartref[]{Gtev}{d} implies
  that
  \begin{equation*}
    \dmapdef{\tev{G\s{F}\s{I}}}{\tp[S]{\Hom{G}{\s{F}}}{\s{I}}}{\Hom{G}{\tp[S]{\s{F}}{\s{I}}}}
  \end{equation*}
  is an isomorphism of $S$--modules. Similar remarks apply to the
  other parts of \thmref[Theorems~]{Ghev} and \thmref[]{Gtev}.
\end{obs}

\begin{ipg}
  Now we turn to applications of \thmref[Theorems~]{Ghev} and
  \thmref[]{Gtev}.
\end{ipg}

\begin{bfhpg}[Remarks]
  The connection between Auslander categories and Gorenstein
  dimensions, as captured by \eqref{main}, allows transfer of
  information between these two realms. This is the theme of section
  \ref{sec:ABstability} and, consequently, each result in that section
  is phrased as either a statement about Auslander categories or a
  statement about Gorenstein dimensions. Thus, the hybrid statements
  in \prpref[]{stabDS} below call for a comment:
  
  Recall that behind the results in \thmref{stabD}, e.g.\ part
  \prtlbl{a}:
  \begin{equation*}
    \Dtp{\F[S]}{\A{D}} \subseteq \A[S]{\s{D}},
  \end{equation*}
  lie stronger inclusions derived from \lemref{HXE}; in this case
  \respartref{SRFC}{a}:
  \begin{equation*}
    \Dtp{\F[S]}{\GF} \subseteq \GF[S].
  \end{equation*}
  Thus, \thmpartref{stabD}{a} could have been phrased as a hybrid,
  \begin{displaymath}
    \tag{\text{$\star$}}
    \Dtp{\F[S]}{\GF} \subseteq \A[S]{\s{D}},
  \end{displaymath}
  and in that form it does not require a dualizing complex for $R$.
  
  In view of \resref{stabcpx} it is natural to seek a result like:
  \begin{displaymath}
    \tag{?}
    \Dtp{\I[S]}{\GF} \subseteq \GI[S] \quad\text{($S$ noetherian).}
  \end{displaymath}
  We do not know if (?) holds in general, but through an application
  of \thmref{Gtev} we obtain the weaker hybrid
  \prppartref[]{stabDS}{a}:
  \begin{equation*}
    \Dtp{\I[S]}{\GF} \subseteq \B[S]{\s{D}}\quad\text{($S$ noetherian
      with a dualizing complex)}.
  \end{equation*}
  Embedded herein are results that can be stated purely in terms of
  Auslander categories or Gorenstein dimensions; we write them out in
  \corref[Corollaries~]{stabDS} and \corref[]{HXF}.

  \begin{spg}
    Note that also $(\star)$, which in the discussion above appears a
    consequence of \eqref{main} and \respartref{SRFC}{a}, can be
    derived easily from \thmref{Gtev}: Let $\s{F} \in \F[S]$ and $A
    \in \GF$.  By \prpcite[(4.4)]{LWC01a} the unit $\unit{\s{F}}$ is
    an isomorphism, and by \resref{stabcpx} the complex
    $\Dtp[S]{\s{D}}{\s{F}}$ belongs to $\I[S]$. Now
    \thmpartref[]{Gtev}{a} and the commutative diagram below shows
    that $\unit{\Dtp{\s{F}}{A}}$ is invertible.
    \begin{displaymath}
      \xymatrix{  \DHom[S]{\s{D}}{\Dtp[S]{\s{D}}{\DtpP{\s{F}}{A}}}
        \ar[rr]^-{\unit{\Dtp{\s{F}}{A}}} & & \Dtp{\s{F}}{A}  \\
        \DHom[S]{\s{D}}{\Dtp{\DtpP[S]{\s{D}}{\s{F}}}{A}}
        \ar[u]_-{\eq}^-{\text{associativity}} 
        \ar[rr]^-{\eq}_-{\tev[SR]{\s{D}\DtpP[S]{\s{D}}{\s{F}}A}} & &
        \Dtp{\DHom[S]{\s{D}}{\Dtp[S]{\s{D}}{\s{F}}}}{A} 
        \ar[u]^-{\eq}_-{\Dtp{\unit{\s{F}}}{A}} 
      }
    \end{displaymath}
    We emphasize that \eqref{main} plays no part in the proof of
    \thmref{Gtev}.
  \end{spg}
\end{bfhpg}

\begin{prp}
  \label{prp:stabDS}
  Let $\fRS$ be a \ho of rings with $\fdf$ finite, and assume that $S$
  is noetherian with a \dc $\s{D}$. The following hold:
  \begin{prt}
  \item If $\s{I}\in\I[S]$ and $A\in\GF$ then $\Dtp{\s{I}}{A} \in
    \B[S]{\s{D}}$
  \end{prt}
  If we also assume that $\F[S] \subseteq \P$, then the following
  hold:
  \begin{prt} \setcounter{prt}{1}
  \item If $A\in\GP$ and $\s{F}\in\F[S]$ then $\DHom{A}{\s{F}} \in
    \A[S]{\s{D}}$
  \item If $\s{I}\in\I[S]$ and $B\in\GI$ then $\DHom{\s{I}}{B} \in
    \A[S]{\s{D}}$
  \end{prt}  
\end{prp}

\begin{prf}
  \proofoftag{c} Under the assumptions we have $\s{I}\in\I$, so
  $\DHom{\s{I}}{B}$ is bounded by \corcite[2.12]{CFH-}. From
  \thmpartref{Ghev}{b} we get an \iso,
  \begin{equation*}
    \Dtp[S]{\s{D}}{\DHom{\s{I}}{B}} \xra[\quad\eq\quad]{\hev{\s{D}\s{I}B}}
    \DHom{\DHom[S]{\s{D}}{\s{I}}}{B}, 
  \end{equation*}
  where the right-hand side is bounded by \corcite[2.12]{CFH-}, as
  \begin{displaymath}
    \DHom[S]{\s{D}}{\s{I}} \in \F[S] \subseteq \P. 
  \end{displaymath}
  by \resref{stabcpx} and the assumptions. Finally, the commutative
  diagram,
  \begin{equation*}
    \xymatrix@C=0.75em{\DHom{\s{I}}{B}
      \ar[d]^-{\eq}_-{\DHom{\counit{\s{I}}}{B}} 
      \ar[rr]^-{\unit{\DHom{\s{I}}{B}}} & &
      \DHom[S]{\s{D}}{\Dtp[S]{\s{D}}{\DHom{\s{I}}{B}}}
      \ar[d]_-{\eq}^-{\DHom[S]{\s{D}}{\hev{\s{D}\s{I}B}}} \\
      \DHom{\Dtp[S]{\s{D}}{\DHom[S]{\s{D}}{\s{I}}}}{B}
      \ar[rr]^-{\eq}_-{\text{adj.}}& & 
      \DHom[S]{\s{D}}{\DHom{\DHom[S]{\s{D}}{\s{I}}}{B}} }
  \end{equation*}
  shows that the unit $\unit{}$ on $\DHom{\s{I}}{B}$ is an \iso in
  $\D[S]$, as desired.
  
  The proofs of the first two parts are similar; part \prtlbl{a}
  relies on \thmpartref{Gtev}{a} and part \prtlbl{b} on
  \thmpartref[]{Gtev}{d}.
\end{prf}

\begin{cor}
  \label{cor:stabDS}
  Let $\fRS$ be a \ho of rings with $\fdf$ finite, and assume that $R$
  and $S$ are noetherian with dualizing complexes $D$ and $\s{D}$,
  respectively. The following hold:
  \begin{prt}
  \item If $\s{I}\in\I[S]$ and $A\in\A{D}$ then $\Dtp{\s{I}}{A} \in
    \B[S]{\s{D}}$
  \item If $A\in\A{D}$ and $\s{F}\in\F[S]$ then $\DHom{A}{\s{F}} \in
    \A[S]{\s{D}}$
  \item If $\s{I}\in\I[S]$ and $B\in\B{D}$ then $\DHom{\s{I}}{B} \in
    \A[S]{\s{D}}$
  \end{prt}  
\end{cor}

\begin{prf}
  Since $R$ is noetherian and has a dualizing complex, we have $\F=\P$
  by \cite[proof of cor.~3.4]{HBF77b} and hence the condition $\F[S]
  \subseteq \P$ is satisfied. The assertions now follow from
  \prpref{stabDS} in view of \eqref{main}.
\end{prf}

\begin{ipg}
  As a second corollary of \prpref[]{stabDS} we get ascent results for
  Gorenstein dimensions:
\end{ipg}

\begin{cor}
  \label{cor:HXF}
  Let $\fRS$ be a \ho of rings with $\fdf$ finite. If $S$ is
  noetherian and has a \dc, then the following hold:
  \begin{prt}
  \item If $\s{I}$ is injective over $S$ and $A$ is \gor flat over
    $R$, then $\tp{\s{I}}{A}$ is \gor injective over $S$
  \end{prt}
  If we also assume that $\F[S] \subseteq \P$, then the following
  hold:
  \begin{prt} \setcounter{prt}{1}
  \item If $A$ is \gor projective over $R$ and $\s{F}$ is flat over
    $S$, then $\Hom{A}{\s{F}}$ is \gor flat over $S$
  \item If $\s{I}$ is injective over $S$ and $B$ is \gor injective
    over $R$, then $\Hom{\s{I}}{B}$ is \gor flat over $S$
  \end{prt}
\end{cor}

\begin{prf}
  The three parts have similar proofs; we write out part \prtlbl{b}:
  By the assumptions $\s{F}\in\P$. For any Gorenstein projective
  $R$--module $A'$, $\DHom{A'}{\s{F}}$ is isomorphic to
  $\Hom{A'}{\s{F}}$ by \corcite[2.10]{CFH-}, so this module has finite
  Gorenstein flat dimension over $S$ by \prppartref{stabDS}{b} and
  \eqref{main}.  Consequently,
  \begin{displaymath}
    \Gfd[S]{\Hom{A'}{\s{F}}}\le d:=\FFD[S]< \infty, 
  \end{displaymath}
  where the first inequality is by \thmcite[3.5]{CFH-} and the second
  follows from \corcite[V.7.2]{rad}. Now consider a piece of a
  complete projective resolution of $A$:
  \begin{equation*}
    0 \to A' \to P_{d-1} \to \cdots \to P_0 \to A \to 0.
  \end{equation*}
  Also $A'$ is Gorenstein projective, and the functor $\Hom{-}{\s{F}}$
  leaves this sequence exact, as $\pd{\s{F}}$ is finite. In the
  ensuing sequence,
  \begin{equation*}
    0 \to \Hom{A}{\s{F}} \to \Hom{P_0}{\s{F}} \to \cdots \to
    \Hom{P_{d-1}}{\s{F}} \to \Hom{A'}{\s{F}} \to 0,
  \end{equation*}
  the $S$--modules $\Hom{P_\l}{\s{F}}$ are flat, cf.~\resref{stab},
  and it follows that $\Hom{A}{\s{F}}$ is Gorenstein flat over $S$ by
  \corcite[3.14]{HHl04a}.
\end{prf}

\vspace{0pt}
\begin{prf}[\bf Proof of \thmref{Ghev}]
  \proofoftag{a} Choose resolutions
  \begin{displaymath}
    \CPb[S] \ni \s{P} \xre Y \quad \text{and} \quad Z \xre \s{J} \in \CIbl[S]
  \end{displaymath}
  and a complex $\CFb[S] \ni \s{F} \eq Z$. Since $R$ is noetherian, we
  can choose a bounded complex $G\eq X$ of finite Gorenstein
  projective $R$--modules. Further, because $\fdf$ is finite, the
  modules in the bounded complex $\Hom[S]{\s{P}}{\s{F}}$ are of finite
  flat dimension over $R$ by \lemref{stabcpxnoeth}.  The modules in
  $G$ are Gorenstein flat by \thmcite[(4.2.6) and (5.1.11)]{LWC}, so
  the complex $\tp{G}{\Hom[S]{\s{P}}{\s{F}}}$ represents
  $\Dtp{X}{\DHom[S]{Y}{Z}}$ by \corcite[2.16]{CFH-}. As $\s{F} \eq
  \s{J}$ there is a \qiso $\s{F} \xre \s{J}$,
  cf.~\rescite[1.4.I]{LLAHBF91}, which by \thmcite[2.15]{CFH-} induces
  another \qiso,
  \begin{displaymath}
    \tag{\one}
    \tp{G}{\Hom[S]{\s{P}}{\s{F}}} \xre \tp{G}{\Hom[S]{\s{P}}{\s{J}}}.
  \end{displaymath}
  Here we use that the modules in $\Hom[S]{\s{P}}{\s{J}}$ have finite
  injective dimension over $R$. Thus, also the right-hand side in
  (\one) represents $\Dtp{X}{\DHom[S]{Y}{Z}}$; furthermore the complex
  $\Hom[S]{\Hom{G}{\s{P}}}{\s{J}}$ represents
  $\DHom[S]{\DHom{X}{Y}}{Z}$ by \corcite[2.10 and rmk.~2.11]{CFH-},
  whence $\hev{XYZ}$ in $\D[S]$ is represented by
  \begin{equation*}
    \tp{G}{\Hom[S]{\s{P}}{\s{J}}} \xra{\hev{G\s{P}\s{J}}}
    \Hom[S]{\Hom{G}{\s{P}}}{\s{J}},
  \end{equation*}
  which is an \iso by \respartref{ev}{b}.
  
  \proofoftag{b} Using that $S$ is noetherian, we choose resolutions
  \begin{displaymath}
    \CfPbr[S] \ni \s{L} \xre U \xre \s{J} \in \CIb[S] 
    \quad \text{and} \quad V \xre \s{I}\in \CIb[S].
  \end{displaymath}
  We also choose a bounded complex $B\eq W$ of Gorenstein injective
  $R$--modules.  The complex $\Hom{\Hom[S]{\s{J}}{\s{I}}}{B}$
  represents $\DHom{\DHom[S]{U}{V}}{W}$ by \corcite[2.12]{CFH-}, as
  all the modules in $\Hom[S]{\s{J}}{\s{I}}$ have finite projective
  dimension over $R$. The last claim relies on \resref{stab} and the
  assumption $\F[S] \subseteq \P$.  By \thmcite[2.9\prtlbl{b}]{CFH-}
  the composite $\s{L} \xre \s{J}$ induces a \qiso,
  \begin{displaymath}
    \Hom{\Hom[S]{\s{L}}{\s{I}}}{B} \xre \Hom{\Hom[S]{\s{J}}{\s{I}}}{B},
  \end{displaymath}
  so also the left-hand side represents $\DHom{\DHom[S]{U}{V}}{W}$.
  Furthermore, $\tp[S]{\s{L}}{\Hom{\s{I}}{B}}$ represents
  $\Dtp[S]{U}{\DHom{V}{W}}$ by \corcite[2.12]{CFH-}, and the \mo
  $\hev[SR]{UVW}$ in $\D[S]$ is represented by
  \begin{equation*}
    \tp[S]{\s{L}}{\Hom{\s{I}}{B}} \xra{\hev[SR]{\s{L}\s{I}B}}
    \Hom{\Hom[S]{\s{L}}{\s{I}}}{B},
  \end{equation*}
  which is an \iso by \respartref{ev}{a}.
\end{prf} 

\vspace{0pt}
\begin{prf}[\vspace*{\thmtopspace} \bf Proof of \thmref{Gtev}]
  \proofoftag{a) and (b} Choose resolutions
  \begin{displaymath}
    \CPbr[S] \ni \s{P} \xre U \xre \s{J}\in\CIb[S] 
    \quad \text{and} \quad V \xre \s{I}\in\CIb[S] 
  \end{displaymath}
  and a bounded complex $A\eq W$ of Gorenstein flat $R$--modules. As
  $\fdf$ is finite and either $S$ or $R$ is noetherian, the bounded
  complex $\Hom[S]{\s{J}}{\s{I}}$ consists of modules of finite flat
  dimension over $R$; cf.~\resref{stab} and \lemref[]{stabcpxnoeth}.
  Therefore, the complex $\tp{\Hom[S]{\s{J}}{\s{I}}}{A}$ represents
  $\Dtp{\DHom[S]{U}{V}}{W}$ by \corcite[2.16]{CFH-}.  The composite
  $\s{P}\xre \s{J}$ induces a \qiso $\Hom[S]{\s{J}}{\s{I}} \xre
  \Hom[S]{\s{P}}{\s{I}}$ between left-bounded complexes of flat and
  injective $S$--modules. By \pgref{hdf} these modules have finite
  flat or finite injective dimension over $R$, so by
  \thmcite[2.15\prtlbl{b}]{CFH-} there is a \qiso
  \begin{displaymath}
    \tp{\Hom[S]{\s{J}}{\s{I}}}{A} \xre \tp{\Hom[S]{\s{P}}{\s{I}}}{A}; 
  \end{displaymath}
  in particular, also the latter complex represents
  $\Dtp{\DHom[S]{U}{V}}{W}$.  Furthermore, \mbox{$\tp{\s{I}}{A}$}
  represents $\Dtp{V}{W}$, again by \corcite[2.16]{CFH-}, so
  $\Hom[S]{\s{P}}{\tp{\s{I}}{A}}$ represents
  $\DHom[S]{U}{\Dtp{V}{W}}$, and the morphism $\tev[SR]{UVW}$ in
  $\D[S]$ is represented by
  \begin{equation*}
    \tp{\Hom[S]{\s{P}}{\s{I}}}{A} \xra{\tev[SR]{\s{P}\s{I}A}}
    \Hom[S]{\s{P}}{\tp{\s{I}}{A}}.
  \end{equation*}
  If $U\in\If[S]$ and $S$ is noetherian, we may assume that
  $\s{P}\in\CfPbr[S]$, and $\tev[SR]{\s{P}\s{I}A}$ is then invertible
  by \respartref{ev}{c}. Similarly, if $R$ is noetherian and
  $W\in\GFf$, we may assume that all modules in $A$ are finite, and
  then $\tev[SR]{\s{P}\s{I}A}$ is invertible by \respartref{ev}{f}.
  
  \proofoftag{c) and (d} Choose complexes
  \begin{displaymath}
    \CFb[S] \ni \s{F} \eq Y 
    \quad \text{ and } \quad
    \CPbr[S] \ni \s{Q} \xre Z \xre \s{J}\in\CIb[S]
  \end{displaymath}
  and a bounded complex $A\eq X$ of Gorenstein projective
  $R$--modules. In the case of \prtlbl{d} we may assume
  $\s{Q}\in\CfPbr[S]$ as $S$ is noetherian; and in the case of
  \prtlbl{c} we may assume that the modules in $A$ are finite as $R$
  is noetherian. We claim that
  \begin{rqm}
  \item the complex $\tp[S]{\Hom{A}{\s{F}}}{\s{Q}}$ represents
    $\Dtp[S]{\DHom{X}{Y}}{Z}$; and
  \item the complex $\Hom{A}{\tp[S]{\s{F}}{\s{Q}}}$ represents
    $\DHom{X}{\Dtp[S]{Y}{Z}}$.
  \end{rqm}
  The modules in $\s{F}$ have finite flat dimension over $R$ by
  \pgref{hdf}, and if $\F[S] \subseteq \P$ they even have finite
  projective dimension over $R$.  Hence (1) follows from \corcite[2.10
  and rmk.~2.11]{CFH-}. To prove (2) we note that as either $S$ or $R$
  is noetherian, \resref{stab} or \lemref[]{stabcpxnoeth} implies that
  all modules in the bounded complex $\tp[S]{\s{F}}{\s{J}}$ have
  finite injective dimension over $R$. Therefore \corcite[2.10]{CFH-}
  gives that the complex $\Hom{A}{\tp[S]{\s{F}}{\s{J}}}$ represents
  $\DHom{X}{\Dtp[S]{Y}{Z}}$.  The composite \qiso $\s{Q}\xre \s{J}$
  induces a \qiso $\tp[S]{\s{F}}{\s{Q}} \xre \tp[S]{\s{F}}{\s{J}}$
  between right-bounded complexes. The modules in
  $\tp[S]{\s{F}}{\s{Q}}$ have finite flat dimension over $R$ according
  to \pgref{hdf}, and if $\F[S] \subseteq \P$ they even have finite
  projective dimension over $R$.  By \thmcite[2.8\prtlbl{b} and
  rmk.~2.11]{CFH-} this \qiso is preserved by $\Hom{A}{-}$.
  
  In total, this shows that the morphism $\tev{XYZ}$ in $\D[S]$ is
  represented by
  \begin{equation*}
    \tp[S]{\Hom{A}{\s{F}}}{\s{Q}} \xra{\tev{A\s{F}\s{Q}}}
    \Hom{A}{\tp[S]{\s{F}}{\s{Q}}}.
  \end{equation*}
  If $\s{Q}\in\CfPbr[S]$ then $\tev{A\s{F}\s{Q}}$ is then invertible
  by \respartref{ev}{e}; this proves \prtlbl{d}. If all the modules in
  $A$ are finite then $\tev{A\s{F}\s{Q}}$ is an \iso by
  \respartref{ev}{d}; this proves \prtlbl{c}.
  
  The proof of \prtlbl{d'} is similar to the proof of \prtlbl{d} and
  thus omitted.
\end{prf}

%%% SECTION 4

\section{Auslander--Buchsbaum formulas}
\label{sec:ABF}

\noindent
We now turn attention to formulas of the Auslander--Buchsbaum type. As
in the previous sections, we consider a relative situation.  That is,
when $X$ is an $R$--complex and $Y$ an $S$--complex, we relate the
$S$--depth of $\Dtp{X}{Y}$ to the depths of $X$ and $Y$ over $R$ and
$S$, respectively.

\begin{bfhpg}[Depth and width]
  The invariants depth and width for complexes over a local ring
  $\Snl$ can be computed in a number of different ways as demonstrated
  in \cite{HBFSIn03}. Here we shall only need two of them.  Let
  $K=\Kc[S]{x_1,\dots,x_e}$ be the Koszul complex on a set of
  generators $x_1,\dots,x_e$ for $\n$. For an $S$--complex $Y$ the
  \emph{depth} and \emph{width} are given by:
  \begin{align}
    \label{eq:dpt}
    \begin{aligned}
      \dpt[S]{Y} &\,=\, -\sup{\DHom[S]{l}{Y}} \\
      &\,=\, -\sup{\DHom[S]{K}{Y}} \\
      &\,=\, e - \sup{\tpP[S]{K}{Y}}, \text{ and}
    \end{aligned}
  \end{align}
  \begin{align}
    \label{eq:wdt}
    \wdt[S]{Y} \,=\, \inf{\DtpP[S]{l}{Y}} \,=\, \inf{\tpP[S]{K}{Y}}.
  \end{align}
  It is easy to prove the inequalities:
  \begin{equation}
    \label{eq:dptineq}
    \dpt[S]{Y} \ge -\sup{Y} \quad \text{ and } \quad \wdt[S]{Y} \ge \inf{Y}.
  \end{equation}
  If $M\ne 0$ is an $\n$--torsion module, that is, \mbox{$M =
    \bigcup_{n=1}^\infty (0:\n^n)_M$}, then $\Hom[S]{l}{M}$ as well as
  $\tp[S]{l}{M}$ are non-zero. This has the following consequence:
  \begin{eqlist}
    \setcounter{equation}{3}
  \item \label{eq:mtor} If $\H{Y}$ is degree-wise $\n$--torsion, then
    there are equalities,
    \begin{displaymath}
      \dpt[S]{Y} = -\sup{Y} \quad \text{ and } \quad \wdt[S]{Y} = \inf{Y}.
    \end{displaymath}
  \end{eqlist}
  For a much stronger statement see \lemcite[2.8]{HBFSIn03}.
\end{bfhpg}

\begin{rmk}
  If $\Rmk$ is local and $X \in \Dfb$, then
  \begin{equation*}
    \tag{\one}
    \pd{X} = \fd{X} = \sup{\DtpP{k}{X}} = - \dpt{\DtpP{k}{X}},
  \end{equation*}
  where the last equality follows from \eqref{mtor}, as all the
  modules in $\H{\Dtp{k}{X}}$ are annihilated by $\m$. The classical
  Auslander--Buchsbaum formula states that if this number (\one) is
  finite, then it equals
  \begin{displaymath}
    \dptR - \dpt{X} = -(\dpt{k} + \dpt{X} - \dptR).
  \end{displaymath}
  Thus one recovers the classical Auslander--Buchsbaum formula by
  setting $Y=k$ and $\f = 1_R$ in \thmref{ABE}\prtlbl{a}. This
  illustrates the point of view that \thmref[]{ABE}\prtlbl{a} is the
  Auslander--Buchsbaum formula for $X$ with coefficients in $Y$ (or
  vice versa).
  
  Given any local homomorphism $\f$, \thmref{ABE} gives
  Auslander--Buchsbaum type formulas for a module $X$ of finite
  classical homological dimension ($\fd[]{}$, $\pd[]{}$, or $\id[]{}$)
  with coefficients in an arbitrary module $Y$.  In \thmref{GABE} we
  relax the conditions on $X$ and, in return, impose conditions on $Y$
  and $\f$ to obtain similar formulas for a modules of finite
  Gorenstein homological dimension.
\end{rmk}

\begin{thm}
  \label{thm:ABE}
  Let $\fRS$ be a local homomorphism of rings. Then the following
  hold:
  \begin{prt}
  \item Let \mbox{$Y \in \Dbl[S]$} and \mbox{$X \in \Dbl$}. If
    \mbox{$\,Y \in \F$} or \mbox{$X \in \F$} then
    \begin{equation*}
      \dpt[S]{\DtpP{Y}{X}} = \dpt[S]{Y} + \dpt{X} - \dptR.
    \end{equation*} 
  \item Let \mbox{$X \in \Dbl$} and \mbox{$Y \in \Dbr[S]$}.  If
    $X\in\P$ or $Y\in\I$ then
    \begin{align*}
      \wdt[S]{\DHom{X}{Y}} = \dpt{X} + \wdt[S]{Y} - \dptR.
    \end{align*}
  \item Let \mbox{$Y \in \Dbl[S]$} and \mbox{$X \in \Dbr$}.  If
    $\,Y\in\P$ or $X\in\I$ then
    \begin{align*}
      \wdt[S]{\DHom{Y}{X}} = \dpt[S]{Y} + \wdt{X} - \dptR.
    \end{align*}
  \end{prt} 

  \begin{spg}
    The absolute version already exists in \thmcite[2.1]{SIn99},
    \cite[(12.8) and (12.20)]{hha}, and \thmcite[(4.14)]{CFF-02}. This
    relative version is joint work between Srikanth Iyengar and Lars
    Winther Christensen. The proof is deferred to the end of this
    section.
  \end{spg}
\end{thm}

\begin{ipg}
  The isomorphisms from \thmref[]{Ghev} and \thmref[]{Gtev} propel the
  main result of this section:
\end{ipg}

\begin{thm}
  \label{thm:GABE}
  Let $\fRS$ be a local homomorphism with $\fdf$ finite.  The
  following hold:
  \begin{prt}
  \item For $Y\in\I[S]$ and $X\in\GF$ there is an equality:
    \begin{align*}
      \dpt[S]{\DtpP{Y}{X}} = \dpt[S]{Y} + \dpt{X} - \dptR.
    \end{align*}        
  \item For $X \in \GP$ and $Y\in\P[S]$ there is an equality:
    \begin{align*}
      \wdt[S]{\DHom{X}{Y}} = \dpt{X} + \wdt[S]{Y} - \dptR.
    \end{align*}
  \item For $Y\in\I[S]$ and $X \in \GI$ there is an equality:
    \begin{align*}
      \wdt[S]{\DHom{Y}{X}} &= \dpt[S]{Y} + \wdt{X} - \dptR.
    \end{align*}
  \end{prt}
\end{thm}

\begin{rmk}
  Comparison of \thmpartref[]{ABE}{a} to \thmpartref[]{GABE}{a} raises
  two questions:
  \begin{itemlist}
  \item Does \thmpartref[]{ABE}{a} hold if one assumes that \mbox{$Y
      \in \F[S]$} instead of \mbox{$Y \in \F$}?
  \item Does \thmpartref[]{GABE}{a} hold without the assumption $Y \in
    \I[S]$?
  \end{itemlist}  
  The answer to the first question is negative: Let $\Rmk$ be
  non-regular and $\f$ be the canonical projection \mbox{$R \onto k$}.
  Then $k \in \F[k]$ but $k \notin \F$, and
  \begin{displaymath}
    \dpt[k]{\DtpP{k}{k}} = -\sup\DHom[k]{k}{\Dtp{k}{k}} = -\sup
    \DtpP{k}{k} = -\infty,
  \end{displaymath}
  while $\dpt[k]{k} + \dpt{k} - \dptR = -\dptR$.
  
  Also the second question has a negative answer: Let $S=R$ be
  Gorenstein but not regular (and $\f$ be the identity map). Then
  \mbox{$k \in \GF$} but \mbox{$k \notin \I$}, and
  \begin{displaymath}
    \dpt{\DtpP{k}{k}} = -\sup{\DtpP{k}{k}} = -\infty,
  \end{displaymath}
  by \eqref{mtor}, while $\dpt{k} + \dpt{k} - \dptR = -\dptR$.
\end{rmk}

\enlargethispage*{0.5cm}
\vspace*{0pt}
\begin{prf}[Proof of \thmref{GABE}]
  \proofoftag{a} We let $\n$ denote the unique maximal ideal of $S$
  and $l$ be the residue field. Furthermore, let $K$ be the Koszul
  complex on a set of generators for $\n$.  In the following sequence
  of isomorphisms, the first and last are by Hom~evaluation
  \respartref{Dev}{a}, and the middle one is by \thmpartref{Gtev}{a}.
  This theorem applies because the complex $\Hom[S]{K}{\E[S]{l}}$ has
  finite injective dimension and finite (length) homology modules:
  \begin{align*}
    \Dtp[S]{K}{\DtpP{\DHom[S]{\E[S]{l}}{Y}}{X}}
    &\eq \Dtp{\DHom[S]{\Hom[S]{K}{\E[S]{l}}}{Y}}{X}\\
    &\eq \DHom[S]{\DHom[S]{K}{\E[S]{l}}}{\Dtp{Y}{X}}\\
    &\eq \Dtp[S]{K}{\DHom[S]{\E[S]{l}}{\Dtp{Y}{X}}}
  \end{align*}
  Because $K$ is depth sensitive, cf.~\eqref{dpt}, this isomorphism
  implies an equality:
  \begin{equation*}
    \tag{\one}
    \dpt[S]{\DtpP{\DHom[S]{\E[S]{l}}{Y}}{X}} \,=\,
    \dpt[S]{\DHom[S]{\E[S]{l}}{\Dtp{Y}{X}}}. 
  \end{equation*}
  The complex $\DHom[S]{\E[S]{l}}{Y}$ is in $\F[S]$, and hence also in
  $\F$ as $\fdf$ is finite. Therefore the left-hand side of (\one) is
  equal to:
  \begin{align*}
    &\dpt[S]{\DHom[S]{\E[S]{l}}{Y}} + \dpt{X} - \dptR \,= \\
    &\wdt[S]{\E[S]{l}} + \dpt[S]{Y} + \dpt{X} - \dptR
  \end{align*}
  by \thmpartref{ABE}{a} and \lemcite[(A.6.4)]{LWC}.  By
  \prpcite[4.6]{HBFSIn03} the right-hand side of (\one) equals
  \begin{align*}
    \wdt[S]{\E[S]{l}} + \dpt[S]{\DtpP{Y}{X}},
  \end{align*}
  and the desired formula follows.
  
  Similar arguments establish parts \prtlbl{b} and \prtlbl{c}: Part
  \prtlbl{b} uses \thmpartref{Gtev}{d'} and \thmpartref[]{ABE}{b},
  while part \prtlbl{c} relies on \thmpartref{Ghev}{b} and
  \thmpartref[]{ABE}{c}.
\end{prf}

\begin{cor}
  \label{cor:SS}
  Let $\fRSl$ be a local homomorphism with $\fdf$ finite.  The
  following hold:
  \begin{prt}
  \item If $X\in\GF$ then $\, \sup{\DtpP{\E[S]{l}}{X}} = \dptR -
    \dpt{X}$.
  \item If $X \in \GP$ then $\, -\inf{\DHom{X}{\Shat}} = \dptR -
    \dpt{X}$.
  \item If $X \in \GI$ then $\, -\inf{\DHom{\E[S]{l}}{X}} = \dptR -
    \wdt{X}$.
  \end{prt}
\end{cor}

\begin{prf}
  For part \prtlbl{a} we set \mbox{$Y=\E[S]{l}$} in
  \thmpartref{GABE}{a} to obtain
  \begin{displaymath}
    -\dpt{\DtpP{\E[S]{l}}{X}} = \dptR - \dpt{X}.
  \end{displaymath}
  Observe that the homology modules of $\Dtp{\E[S]{l}}{X}$ are
  $\m$--torsion, as $\f$ is local and $\E{l}$ is $\n$--torsion, and
  apply \eqref{mtor}.  Part \prtlbl{b} follows from \prtlbl{a} as
  \mbox{$\GP \subseteq \GF$}:
  \begin{align*}
    -\inf{\DHom{X}{\Shat}}
    &= -\inf{\DHom{X}{\Hom[S]{\E[S]{l}}{\E[S]{l}}}} \\
    &= -\inf{\DHom[S]{\Dtp{X}{\E[S]{l}}}{\E[S]{l}}} \\
    &= \sup{\DtpP{\E[S]{l}}{X}}.
  \end{align*}
  To prove \prtlbl{c} we need the representations of local
  (co)homology from \cite{AJL-97,JPGJPM92} as summed up in
  \rescite[(2.6)]{AFr03}:
  \begin{displaymath}
    \DL[\n]{\DHom{\E[S]{l}}{X}} \eq \DHom{\DG[\n]{\E[S]{l}}}{X} \eq
    \DHom{\E[S]{l}}{X}.
  \end{displaymath}
  Consequently,
  \begin{align*}
    -\inf{\DHom{\E[S]{l}}{X}} &= -\infP{\DL[\n]{\DHom{\E[S]{l}}{X}}} \\
    &=  -\wdt[S]{\DHom{\E[S]{l}}{X}} \\
    &= \dptR - \wdt{X}
  \end{align*}
  by \thmcite[(2.11)]{AFr03} and \thmpartref{GABE}{c}.
\end{prf}

\begin{ipg}
  A finite Gorenstein projective module $G$ over a local ring $R$ has
  $\dpt{G} = \dptR$. The next corollary extends this equality to
  non-finite modules.
\end{ipg}

\begin{cor}
  Let $\Rmk$ be local. If $A$ is a Gorenstein flat module of finite
  depth, then $\dpt{A} =  \dptR$.  Similarly, if $B$ is a \gor
  injective modules of finite width, then $\wdt{B} =  \dptR$.
\end{cor}

\begin{prf}
  By \corpartref[]{SS}{c} we have $-\inf{\DHom{\E{k}}{B}} = \dptR -
  \wdt{B}$. If $\wdt{B}$ is finite, so is $-\inf{\DHom{\E{k}}{B}}$,
  and since $\DHom{\E{k}}{B}$ is represented by the module
  $\Hom{\E{k}}{B}$, cf.~\corcite[2.12]{CFH-}, the infimum must be
  zero. This proves the second statement; the proof of the first one
  is similar.
\end{prf}

\enlargethispage*{0.7cm}
\begin{ipg}
  We close this section with the proof of \thmref{ABE}. It is broken
  down into six steps, the first of which is a straightforward
  generalization of the argument in \thmcite[2.4]{HBFSIn03} to the
  relative situation. The third step uses \respartref{Dev}{f} and is
  considerably shorter than the proof of the absolute version in
  \thmcite[(4.14)]{CFF-02}.
\end{ipg}

\begin{prf}[Proof of \thmref{ABE}]
  \step{1} First we assume that \mbox{$X\in\F$} and \mbox{$Y \in
    \Dbl[S]$}. The second equality in the next computation follows by
  tensor evaluation \respartref{Dev}{d}; the third holds as all
  $l$--vector spaces become vector spaces over $k$ through the local
  \ho $\f$,

  \begin{align*}
    \dpt[S]{\DtpP{Y}{X}} &= -\sup{\DHom[S]{l}{\Dtp{Y}{X}}} \\
    &= -\supP{\Dtp{\DHom[S]{l}{Y}}{X}} \\
    &= -\supP{\tp[k]{\DHom[S]{l}{Y}}{\DtpP{k}{X}}} \\
    &= -\sup \DHom[S]{l}{Y} - \sup{\DtpP{k}{X}} \\
    &= \dpt[S]{Y} - \sup{\DtpP{k}{X}}.
  \end{align*}
  In particular, for $S=R=Y$ and $\f=1_R$ we get $\dpt{X} = \dptR -
  \sup{\DtpP{k}{X}}$; combining this with the equality above, the
  desired equality follows.

  \enlargethispage*{0.5cm}
  \step{2} Next we assume that \mbox{$Y\in\F$} and \mbox{$X \in
    \Dbl$}. Let $K = \Kc[S]{x_1,\dots,x_e}$ be the Koszul complex on a
  set of generators for $\n$. In the next computation, the first and
  last equalities are by \eqref{dpt}, while the second and penultimate
  ones follow by \eqref{mtor}. (Since $\f$ is local, the homology
  modules of $\tp[S]{K}{Y}$ and $\Dtp[S]{K}{\DtpP{Y}{X}} \eq
    \Dtp{\tpP[S]{K}{Y}}{X}$ are $\m$--torsion, even annihilated by
  $\m$ cf.~\thmcite[16.4]{Mat}.)  The Koszul complex $K$ consists of
  finite free $S$--modules, so also $\tp[S]{K}{Y}$ belongs to $\F$;
  the third equality below is therefore the absolute version of the
  formula already established in \step{1}.
  \begin{align*}
    \dpt[S]{\DtpP{Y}{X}} &= e - \sup{\DtpP[S]{K}{\DtpP{Y}{X}}} \\
    &= e + \dpt{\DtpP{\tpP[S]{K}{Y}}{X}} \\
    &= e + \dpt{\tpP[S]{K}{Y}} + \dpt{X} - \dptR \\
    &= e - \sup{\tpP[S]{K}{Y}} + \dpt{X} - \dptR \\
    &= \dpt[S]{Y} + \dpt{X} - \dptR.
  \end{align*}
  This concludes the proof of part \prtlbl{a}. The arguments
  establishing \prtlbl{b} and \prtlbl{c} are intertwined, and the
  whole argument is divided into four steps (\step{3}--\step{6}).
  
  \step{3} We establish the equality in \prtlbl{b} under the
  assumption that $X\in\P$ and $Y \in \Dbr[S]$. The second equality in
  the next computation follows by tensor evaluation,
  \respartref{Dev}{f}, and the third by adjunction:
  \begin{align*}
    \wdt[S]{\DHom{X}{Y}} &= \infP{\Dtp[S]{\DHom{X}{Y}}{l}} \\
    &= \inf \DHom{X}{\Dtp[S]{Y}{l}} \\
    &= \inf \Hom[k]{\Dtp{X}{k}}{\Dtp[S]{Y}{l}} \\
    &= \inf{\DtpP[S]{Y}{l}} - \sup{\DtpP{X}{k}}.
  \end{align*}
  Since $X\in\P \subseteq \F$ we have $-\sup{\DtpP{X}{k}} = \dpt{X} -
  \dptR$, as established in \step{1} above, and by \eqref{wdt} we have
  $\inf{\DtpP[S]{Y}{l}} = \wdt[S]{Y}$.
  
  \step{4} Based on \step{3} we can prove the equality in \prtlbl{c}
  under the assumptions that $Y\in\P$ and $X \in \Dbr$:
  \begin{align*}
    \wdt[S]{\DHom{Y}{X}} &= \inf{\DtpP[S]{K}{\DHom{Y}{X}}} \\
    &= \wdt{\DtpP[S]{K}{\DHom{Y}{X}}} \\
    &= \wdt{\DHom{\Hom[S]{K}{Y}}{X}} \\
    &= \dpt{\Hom[S]{K}{Y}} + \wdt{X} - \dptR \\
    &= -\sup{\Hom[S]{K}{Y}} + \wdt{X} - \dptR \\
    &= \dpt[S]{Y} + \wdt{X} - \dptR.
  \end{align*}

  The first equality is by \eqref{wdt}, the second by \eqref{mtor},
  and the third by Hom~evaluation \respartref{Dev}{a}. The complex
  $\Hom[S]{K}{Y}$ is in $\P$, so the fourth equality follows from what
  we have already proved in \step{3} with \mbox{$\f=1_R$}. The
  penultimate equality is by \eqref{mtor}, as $\Hom[S]{K}{Y}$ is
  shift-isomorphic to $\Dtp[S]{K}{Y}$, and hence its homology is
  annihilated by $\m$. The last equality is by \eqref{dpt}.
  
  \step{5} We prove \prtlbl{c} under the assumption that \mbox{$X \in
    \I$} and \mbox{$Y \in \Dbl[S]$}.  The second equality below
  follows by Hom~evaluation \respartref{Dev}{b} and the third by
  adjunction,
  \begin{align*}
    \wdt[S]{\DHom{Y}{X}} &= \inf{\DtpP[S]{l}{\DHom{Y}{X}}} \\
    &= \inf{\DHom{\DHom[S]{l}{Y}}{X}} \\
    &= \inf{\Hom[k]{\DHom[S]{l}{Y}}{\DHom{k}{X}}} \\
    &= \inf{\DHom{k}{X}} - \sup{\DHom[S]{l}{Y}} \\
    &= \inf{\DHom{k}{X}} + \dpt[S]{Y}.
  \end{align*}
  In particular, with $S=R=Y$ and $\f=1_R$ we get $\wdt{X} =
  \inf{\DHom{k}{X}} + \dptR$, and the desired formula follows. This
  concludes to proof of part \prtlbl{c}.
  
  \step{6} A computation similar to the one performed in \step{4} ---
  but this time based on \respartref{Dev}{e} and the absolute case
  $\f=1_R$ of \step{5} --- proves \prtlbl{b} under the assumptions
  that $Y \in \I$ and $X \in \Dbl$. This concludes the proof.
\end{prf}

%%% SECTION 5

\section{Catalogues}
\label{sec:cat}

\noindent
In this final section we catalogue ascent properties of Auslander
categories and Gorenstein dimensions. It summarizes the results proved
in sections \ref{sec:ABstability} and \ref{sec:applications} and fills
the gaps that become apparent when the results are presented
systematically.

\begin{bfhpg}[Ascent cross tables]
  \renewcommand{\arraystretch}{0.85} Let $\fRS$ be a \ho of rings. The
  cross tables below sum up ascent properties of Auslander categories.
  
  The results from sections \ref{sec:ABstability} and
  \ref{sec:applications} gives general results for half of the
  combinations considered in these tables. For the other half, we
  provide references to counterexamples (cntrex) and in some cases to
  interesting special cases (sp case).
  
  \begin{spg}
    In \pgref{tab:T} we assume that $S$ is noetherian and $\s{C}$ a
    \sdc for $S$.

    \begin{stabmatrix}{T}
      \tag{a}
      \begin{tabular}{r|cc}
        \sym{\Dtp{-}{-}} & \sym{\F} & \sym{\I} \\[0.5ex]

        \hline \\[-1ex]

        \sym{\A[S]{\s{C}}} & \sym{\A[S]{\s{C}}} & cntrex/sp case \\

        \cit{$S$ noetherian} & \cit{by \prppartref[]{ABstabC}{a}} &
        \cit{see \obsref[]{gor}/\obspartref[]{gor}{a}} \\[2ex]

        \sym{\B[S]{\s{C}}} & \sym{\B[S]{\s{C}}} & \text{cntrex} \\

        \cit{$S$ noetherian} & \cit{by \prppartref[]{ABstabC}{b}} &
        \cit{see \exapartref[]{counter}{a}}
      \end{tabular}     
    \end{stabmatrix}

    \noindent In \pgref{tab:Tf} we assume that $R$ and $S$ are
    noetherian, $\fdf$ is finite, and $D$ and $\s{D}$ are dualizing
    complexes for $R$ and $S$, respectively.

    \begin{stabmatrix}{Tf}
      \tag{b}
      \begin{tabular}{r|cc}
        \sym{\Dtp{-}{-}} & \sym{\A{D}} & \sym{\B{D}} \\  

        \cit{$\fdf$ finite}\hspace{0.05ex} & \cit{$R$ noetherian} & \cit{$R$
          noetherian} \\[0.5ex] 

        \hline \\[-1ex]

        \sym{\F[S]} & \sym{\A[S]{\s{D}}} & cntrex/sp case \\  

        \cit{$S$ noetherian} & \cit{by \thmpartref[]{stabD}{a}} &
        \cit{see \exapartref[]{countermap}{a}/\obspartref[]{gor}{b}} \\[2ex] 

        \sym{\I[S]} & \sym{\B[S]{\s{D}}} & \text{cntrex} \\

        \cit{$S$ noetherian} & \cit{by \corpartref[]{stabDS}{a}} &
        \cit{see \exapartref[]{counter}{a}}
      \end{tabular}
    \end{stabmatrix}

    \noindent
    In \pgref{tab:H2} and \pgref{tab:H1} we assume that $S$ is
    noetherian and $\s{C}$ semi-dualizing for $S$.

    \begin{stabmatrix}{H2}
      \tag{c}
      \begin{tabular}{r|cc}
        \sym{\DHom{-}{-}} & \sym{\A[S]{\s{C}}} &  \sym{\B[S]{\s{C}}} \\  

        {} & \cit{$S$ noetherian} & \cit{$S$ noetherian} \\[0.5ex] 

        \hline \\[-1ex]

        \sym{\P} & \sym{\A[S]{\s{C}}} & \sym{\B[S]{\s{C}}} \\

        {} & \cit{by \prppartref[]{ABstabC}{c}} & \cit{by
          \prppartref[]{ABstabC}{d}} \\[2ex] 

        \sym{\I} & \text{cntrex} & cntrex/sp case \\  

        {} & \cit{see \exapartref[]{counter}{b}} & \cit{see
          \obsref[]{gor}/\obspartref[]{gor}{a}}
      \end{tabular}
    \end{stabmatrix}

    \begin{stabmatrix}{H1}
      \tag{d}
      \begin{tabular}{r|cc}
        \sym{\DHom{-}{-}} & \sym{\P} &  \sym{\I} \\[0.5ex]

        \hline \\[-1ex]

        \sym{\A[S]{\s{C}}} & cntrex/sp case & \sym{\B[S]{\s{C}}} \\

        \cit{$S$ noetherian}& \cit{see \obsref[]{gor}/\obspartref[]{gor}{a}} & \cit{by
          \prppartref[]{ABstabC}{e}} \\[2ex]

        \sym{\B[S]{\s{C}}} & \text{cntrex} & \sym{\A[S]{\s{C}}} \\

        \cit{$S$ noetherian}& \cit{see \exapartref[]{counter}{b}} &
        \cit{by \prppartref[]{ABstabC}{f}}
      \end{tabular}
    \end{stabmatrix}
    
    \noindent
    In \pgref{tab:H2f} and \pgref{tab:H1f} we assume that $R$ and $S$
    are noetherian, $\fdf$ is finite, and $D$ and $\s{D}$ are dualizing
    complexes for $R$ and $S$, respectively.

    \begin{stabmatrix}{H2f}
      \tag{e}
      \begin{tabular}{r|cc}
        \sym{\DHom{-}{-}} & \sym{\A{D}} & \sym{\B{D}} \\  

        \cit{$\fdf$ finite}\hspace{2em} & \cit{$R$ noetherian} & \cit{$R$
          noetherian} \\[0.5ex] 

        \hline \\[-1ex]

        \sym{\P[S]} & cntrex/sp case & \sym{\B[S]{\s{D}}} \\

        \cit{$S$ noetherian} & \cit{see
          \exapartref[]{countermap}{b}/\obspartref[]{gor}{b}} & \cit{by
          \thmpartref[]{stabD}{b}} \\[2ex] 

        \sym{\I[S]} & \text{cntrex} & \sym{\A[S]{\s{D}}} \\

        \cit{$S$ noetherian} & \cit{see \exapartref[]{counter}{b}} &
        \cit{by \corpartref[]{stabDS}{c}}   
      \end{tabular}
    \end{stabmatrix}

    \begin{stabmatrix}{H1f}
      \tag{f}
      \begin{tabular}{r|cc}
        \sym{\DHom{-}{-}} & \sym{\P[S]} &  \sym{\I[S]} \\

        \cit{$\fdf$ finite}\hspace{2em} & \cit{$S$ noetherian} & \cit{$S$
          noetherian} \\[0.5ex]

        \hline \\[-1ex]

        \sym{\A{D}} & \sym{\A[S]{\s{D}}} &  \sym{\B[S]{\s{D}}} \\

        \cit{$R$ noetherian} & \cit{by \corpartref[]{stabDS}{b}} & \cit{by
          \thmpartref[]{stabD}{c}} \\[2ex]

        \sym{\B{D}} & \text{cntrex} & cntrex/sp case \\

        \cit{$R$ noetherian} & \cit{see \exapartref[]{counter}{b}} &
        \cit{see \exapartref[]{countermap}{c}/\obspartref[]{gor}{b}} 
      \end{tabular}
    \end{stabmatrix}
  \end{spg}
\end{bfhpg}

\begin{obs}
  \label{obs:gor}
  Let $R$ be noetherian and $C$ be a \sdc for $R$.  In general, the
  combinations
  \begin{equation*}
    \Dtp{\A{C}}{\I},\quad \DHom{\A{C}}{\P}\quad \text{and}\quad
    \DHom{\I}{\B{C}}
  \end{equation*}
  do not even end up in $\Db$, and in particular not in $\A{C}$ or
  $\B{C}$.  Indeed, let $R$ be the ring from \exaref{counter} and set
  $C=R$. The \dm $D$ belongs to both $\I$ and $\A{R} = \Db$, and $R$
  belongs to both $\P$ and $\B{R} = \Db$, but $\Dtp{D}{D}$ as well as
  $\DHom{D}{R}$ is unbounded.
  
  \begin{spg}
    \prtlbl{a}: Let $\fRS$ be a homomorphism of noetherian rings. Let
    $D$ be a dualizing complex for $R$ and assume that
    $\s{D}=\Dtp{S}{D}$ is dualizing\footnote{Local
      homomorphisms with this property are called
      \emph{quasi-Gorenstein} and were first studied in
      \cite{LLAHBF97}.} for $S$.  An $S$--complex is then in $\A[S]{\s{D}}$ or
    $\B[S]{\s{D}}$ if and only if it belongs to $\A{D}$ or $\B{D}$,
    respectively; see \cite[proof of cor.~(7.9)]{LLAHBF97}.  Applying
    \corref{stabDS} to $\f=1_R$ now yields
    \begin{itemlist}
    \item $\Dtp{\A[S]{\s{D}}}{\I} \subseteq \B[S]{\s{D}}$,
    \item $\DHom{\A[S]{\s{D}}}{\P} \subseteq \A[S]{\s{D}}$, and
    \item $\DHom{\I}{\B[S]{\s{D}}} \subseteq \A[S]{\s{D}}$.
    \end{itemlist}
  \end{spg}
  
  \begin{spg}
    \prtlbl{b}: Let $\fRS$ be a homomorphism of noetherian rings with
    $\fdf$ finite. If $R$ has a dualizing complex, $D$, then $\P=\F$,
    e.g.\ by \cite[proof of cor.~3.4]{HBF77b}, and hence it is
    immediate by the absolute version of
    \prppartref[Prop.~]{ABstabC}{b,c,f}, cf.~\pgref{hdf}, that
    \begin{itemlist}
    \item $\Dtp{\F[S]}{\B{D}} \subseteq \B{D}$,
    \item $\DHom{\P[S]}{\A{D}} \subseteq \A{D}$, and
    \item $\DHom{\B{D}}{\I[S]} \subseteq \A{D}$.
    \end{itemlist}
    Under the additional assumption that $\s{D} = \Dtp{S}{D}$ is
    dualizing\footnote{Local homomorphisms with this property are
      called \emph{Gorenstein} and were first studied in
      \cite{LLAHBF90}.}  for $S$, it follows, as above, that e.g.\ the
    combination $\Dtp{\F[S]}{\B{D}}$ even ends up in
    $\B[S]{\s{D}}$.
  \end{spg}
  \noindent It is easy to see that this is not the general behavior:
\end{obs}

\begin{exa}
  \label{exa:countermap}
  Let $R$ be a field and consider the flat map
  \begin{equation*}
    \dmapdef{\f}{R}{S=\pows[R]{X,Y}/(X^2\!,XY,Y^2)}.
  \end{equation*}
  $R$ is Gorenstein, so $D=R$ is dualizing for $R$ and
  $\A{D}=\B{D}=\Db$. The ring $S$ has a dualizing complex $\s{D} \eq
  \DHom{S}{R}$ but is not Gorenstein, see \exaref{counter}. In
  particular, $S$ is not in $\B[S]{\s{D}}$, and $\s{D}$ is not in
  $\A[S]{\s{D}}$. Nevertheless,
  \begin{align*}
    \tag{a} \Dtp{\F[S]}{\B{D}} \ni&\; \Dtp{S}{R} \eq S \not\in \B[S]{\s{D}},\\
    \tag{b} \DHom{\P[S]}{\A{D}} \ni&\: \DHom{S}{R} \not\in \A[S]{\s{D}},\ \text{and}\\
    \tag{c} \DHom{\B{D}}{\I[S]} \ni&\: \DHom{R}{\s{D}} \eq \s{D}
    \not\in \A[S]{\s{D}}.
  \end{align*}
\end{exa}

\begin{rmk}
  In general one cannot expect ascent results involving two modules of
  finite Gorenstein dimension or two Auslander categories: Let $R$ be
  a local non-regular \gor ring with residue field $k$, e.g.\ 
  $R=\poly{X}/(X^2)$. In this case, $k$ has finite Gorenstein
  projective dimension and finite Gorenstein injective dimension; in
  particular, $k\in\A{D} = \B{D} = \Db$, but both $\Dtp{k}{k}$ and
  $\DHom{k}{k}$ are unbounded as $\pd{k} = \fd{k} = \id{k} = \infty$.
\end{rmk}
\clearpage
\newlength{\miniindent}
\setlength{\miniindent}{\parindent}
\begin{sideways}
  \begin{minipage}{\textheight}
    \begin{bfhpg}[Ascent table I]
      \renewcommand{\arraystretch}{1.25}
      \label{table1}
      Let $\fRS$ be a \ho of rings. The table below sums up ascent
      results for Gorenstein dimensions pertaining to that situation. We
      use abbreviated notation, e.g.\ ``G--projective/S'' means a
      Gorenstein projective $S$--module. Note that the ``G'' always
      lives over $S$.

      \begin{longtable}[c]{crcl@{\hspace{2em}}l}
        & Ascent result & & &\hspace{6em} Requirements \\
        \hline
        \TAF & \Ttp{\GFS}{\FR} & is & \sl \GFS \\
        \TXP & \Ttp{\GPS}{\PR} & is & \sl \GPS \\
        \TBF & \Ttp{\GIS}{\FR} & is & \sl \GIS & $S$ is noetherian and
        has a \dc \\
        \TBL & \Ttp{\GIS}{\pR} & is & \sl \GIS \\[1.5ex]
        
        \HFBB & \THom{\FR}{\GIS} & is & \sl \GIS\ & $\F = \P$ or $\F[S]
        = \P[S]$ \\
        \HPB & \THom{\PR}{\GIS} & is & \sl \GIS \\
        \HPA & \THom{\PR}{\GFS} & is & \sl \GFS & $S$ is noetherian and
        has a \dc \\
        \HLA & \THom{\pR}{\GFS} & is & \sl \GFS \\
        \HLX & \THom{\pR}{\GPS} & is & \sl \GPS \\[1.5ex]

        \HAE & \THom{\GFS}{\IR} & is & \sl \GIS \\
        \HBE & \THom{\GIS}{\IR} & is & \sl \GFS & $S$ is noetherian and
        has a \dc \\
      \end{longtable}
    \end{bfhpg}

    \begin{spg}
      {\bf Comments:} The requirements in \HFBB\ hold if $R$ or $S$ is
      noetherian with finite Krull dimension; see footnote \ref{fnt}
      to \thmref{Ghev}.
      
      \hspace{\miniindent} In the absolute case, $\f=1_R$, the two
      requirements in \pgref{table1}\HFBB\ are the same.
    \end{spg}
  \end{minipage}
\end{sideways}  

\clearpage

\begin{sideways}  
  \begin{minipage}{\textheight}
    \begin{bfhpg}[Ascent table II]
      \renewcommand{\arraystretch}{1.25}
      \label{table2}
      Let $\fRS$ be a \ho with $\fdf$ finite. The table below sums up
      ascent results pertaining to that situation.  We use the
      abbreviated notation from \pgref{table1}; note that in this table
      the ``G'' always lives over $R$.

      \begin{longtable}[c]{crcl@{\hspace{2em}}l}
        & Ascent result & & & \hspace{4em} Requirements \\
        \hline
        \TFA & \Ttp{\FS}{\GFR} & is & \sl \GFS \\
        \TFX & \Ttp{\FS}{\GPR} & is & \sl \GFS & $S$ is noetherian and
        $\F[S] \subseteq \P$ \\
        \TPX & \Ttp{\PS}{\GPR} & is & \sl \GPS & $\P[S] \subseteq \P$ \\
        \TEA & \Ttp{\IS}{\GFR} & is & \sl \GIS & $S$ is noetherian with
        a \dc \\
        \TEG & \Ttp{\IS}{\GpR} & is & \sl \GIS & $R$ is noetherian
        \\[1.5ex]
        \HAEE & \THom{\GFR}{\IS} & is & \sl \GIS \\
        \HXEE & \THom{\GPR}{\IS} & is & \sl \GIS & $S$ is noetherian and
        $\F[S] \subseteq \P$ \\
        \HGE & \THom{\GpR}{\IS} & is & \sl \GIS & $R$ is noetherian \\
        \HXF & \THom{\GPR}{\FS} & is & \sl \GFS & $S$ is noetherian with
        a \dc and $\F[S] \subseteq \P$ \\
        \HGF & \THom{\GpR}{\FS} & is & \sl \GFS & $R$ is noetherian \\
        \HGP & \THom{\GpR}{\PS} & is & \sl \GPS & $R$ is noetherian
        \\[1.5ex]
        \HFB & \THom{\FS}{\GIR} & is & \sl \GIS & (1) $\F[S] = \P[S]
        \subseteq \P$; or \\
        & & & & (2)  $\F[S] \subseteq \P$, and $R$ or $S$ is noetherian \\
        \HPBB & \THom{\PS}{\GIR} & is & \sl \GIS & $\P[S] \subseteq \P$ \\
        \HEB & \THom{\IS}{\GIR} & is & \sl \GFS & $S$ is noetherian with
        a \dc and $\F[S] \subseteq \P$ \\
      \end{longtable}
    \end{bfhpg}
    \begin{spg}
      {\bf Comments:} Part \TPX\ is an unpublished result of
      Hans-Bj{\o}rn Foxby.
      
      \hspace{\miniindent}Note that with the exception of \TEA, \HXF\ 
      and \HEB, which require a \dc for $S$, all results in this table
      hold when $R$ and $S$ are noetherian and $R$ has finite Krull
      dimension; see the footnote \ref{fnt} to \thmref{Ghev}.
      
      \hspace{\miniindent}In the absolute case, $\f=1_R$, the first
      requirement in \pgref{table2}\HFB\ is weaker than the second.
    \end{spg}
  \end{minipage}
\end{sideways}  

\begin{bfhpg}[Remarks]
  Four results in Table II, namely \pgpartref{table2}{c',e',f',g},
  deal with finite \gor projective modules over a noetherian ring.
  Recall that a finite module over such a ring is \gor projective if
  and only if it is \gor flat if and only if it is totally reflexive
  (in the sense of \cite{LLAAMr02}), cf.~\thmcite[(4.2.6) and
  (5.1.11)]{LWC}.
  
  \begin{spg}
    Over a noetherian ring where every flat module has finite
    projective dimension, e.g.\ a ring of finite Krull dimension,
    every \gor projective module is \gor flat,
    cf.~\prpcite[3.4]{HHl04a}.  Thus, for example, \pgref{table1}\TAF\ 
    includes the result:
    \begin{center}\sl \small
      \Ttp{\GPS}{\FR}\ is \textsl{\GFS} when $S$ is noetherian and
      $\F[S] = \P[S]$.
    \end{center}
    By the same token, \pgref{table2}\TFX\ also holds when $R$ is
    noetherian with $\F=\P$; that setting it is a special case of
    \pgref{table2}\TFA.
  \end{spg}
\end{bfhpg}

\vspace*{0pt}
\begin{prf}[\bf Proofs for ascent table I]
  Parts \TBF, \HPA, and \HBE\ are proved in \prpref[Prop.~]{HPA}.
  
  The proofs of parts \TAF, \TXP, \TBL, \HPB, \HLA, \HLX, and \HAE\ 
  are entirely functorial and follow the pattern from the proof of
  \lemref{HXE}. 
  
  This leaves only \HFBB: Let $F$ be flat over $R$ and $\s{B}$ be \gor
  injective over $S$. Under either assumption, the module $\tp{F}{S}$
  has finite projective dimension over $S$, so there exists an exact
  sequence $0 \to \s{P}_d \to \cdots \to \s{P}_0 \to \tp{F}{S} \to 0$,
  where the $\s{P}$'s are projective over $S$.  The functor
  $\Hom[S]{-}{\s{B}}$ leaves this sequence exact,
  \begin{align*}
    0 \to \Hom[S]{\tp{F}{S}}{\s{B}} \to \Hom[S]{\s{P}_0}{\s{B}} \to
    \cdots \to \Hom[S]{\s{P}_d}{\s{B}} \to 0.
  \end{align*}
  It follows that \mbox{$\Gid[S]{\Hom[S]{\tp{F}{S}}{\s{B}}}\le d$}, as
  each module $\Hom[S]{\s{P}_\l}{\s{B}}$ is Gorenstein injective over
  $S$ by \HPB; note that $d$ is independent of $\s{B}$.  Next,
  consider a piece of a complete injective resolution of $\s{B}$, say
  $0 \to \s{B}'\to \s{I}_{d-1} \to \cdots \to \s{I}_0 \to \s{B} \to
  0$. Also $\s{B}'$ is \gor injective, so
  \mbox{$\Gid[S]{\Hom[S]{\tp{F}{S}}{\s{B}'}} \le d$}, and applying
  $\Hom[S]{\tp{F}{S}}{-}$ we get the exact sequence
  \begin{multline*}
    0 \to \Hom[S]{\tp{F}{S}}{\s{B}'} \to
    \Hom[S]{\tp{F}{S}}{\s{I}_{d-1}}
    \to \cdots\\
    \cdots \to \Hom[S]{\tp{F}{S}}{\s{I}_0} \to
    \Hom[S]{\tp{F}{S}}{\s{B}} \to 0.
  \end{multline*}
  The modules $\Hom[S]{\tp{F}{S}}{\s{I}_{\l}}$ are injective over $S$,
  by \resref{stab}, and hence $\Hom[S]{\tp{F}{S}}{\s{B}} \is
  \Hom{F}{\s{B}}$ is \gor injective over $S$.
\end{prf}

\vspace*{0pt}
\begin{prf}[\bf Proofs for ascent table II]
  Parts \TEA, \HXF, and \HEB\ are proved in \corref{HXF}; parts \TFA,
   \HAEE, and \HPBB\ are proved in \lemref{HXE}.
  
  \TFX: Let $\crs{P}$ be complete projective resolution over $R$ and
  $\s{F}$ a flat $S$--module. The complex $\tp{\s{F}}{\crs{P}}$ is
  exact and consists of flat $S$--modules.  Now let $\s{J}$ be any
  $S$--injective module; by adjunction we have
  \begin{align*}
    \Hom[\ZZ]{\tp[S]{\s{J}}{\tpP{\s{F}}{\crs{P}}}}{\QQ/\ZZ} \is
    \Hom{\crs{P}}{\Hom[\ZZ]{\tp[S]{\s{J}}{\s{F}}}{\QQ/\ZZ}}.
  \end{align*}
  Since $S$ is noetherian, the module
  $\Hom[\ZZ]{\tp[S]{\s{J}}{\s{F}}}{\QQ/\ZZ}$ is $S$--flat and hence in
  $\P$.  This shows that $\tp{\s{F}}{\crs{P}}$ is a complete flat
  resolution over $S$.
  
  \TPX: Let $\s{P}$ be a projective $S$--module and $\crs{P}$ a
  complete projective resolution over $R$.  The complex
  $\tp{\s{P}}{\crs{P}}$ consists of projective $S$--modules, and it is
  exact, as $\s{P}\in\P[S] \subseteq \F$. For any projective
  $S$--module $\s{Q}$, we have $\Hom[S]{\tp{\s{P}}{\crs{P}}}{\s{Q}}
  \is \Hom[S]{\s{P}}{\Hom{\crs{P}}{\s{Q}}}$, which is exact as
  $\s{Q}\in\P[S] \subseteq \P$.
  
  \TEG: Let $\crs{L}$ be a complete resolution by finite free
  $R$--modules ($R$ is noetherian). By \lemcite[(5.1.10)]{LWC},
  $\crs{L}$ is also a complete flat resolution.  Since $\fdf$ is
  finite, the injective dimension of $\s{I}$ over $R$ is finite, so
  $\tp{\s{I}}{\crs{L}}$ is an exact complex of injective $S$--modules.
  To see that $\tp{\s{I}}{\crs{L}}$ is a complete injective resolution
  over $S$, let $\s{J}$ be any injective $S$--module, and apply
  \respartref{ev}{e} to obtain
  \begin{equation*}
    \tag{\one}
    \Hom[S]{\s{J}}{\tp{\s{I}}{\crs{L}}} \,\is\,
    \tp{\Hom[S]{\s{J}}{\s{I}}}{\crs{L}}. 
  \end{equation*}
  By \lemref{stabcpxnoeth} the module $\Hom[S]{\s{J}}{\s{I}} \eq
  \DHom[S]{\s{J}}{\s{I}}$ has finite flat dimension over $R$. The
  complete flat resolution $\crs{L}$ remains exact when tensored by a
  module in $\F$, cf.~\lemcite[2.3]{CFH-}, so the complex in (\one) is
  exact.
  
  \HXEE: Let $G$ be a \gor projective $R$--module and $\s{I}$ an
  injective $S$--module. Under the assumptions, $\tp{G}{S}$ is \gor
  flat over $S$ by \TFX\ and $\Hom{G}{\s{I}} \is
  \Hom[S]{\tp{G}{S}}{\s{I}}$ is then \gor injective over $S$ by
  \pgref{table1}\HAE.
  
  \HGE: Let $\crs{L}$ be a complete resolution by finite free
  $R$--modules and $\s{I}$ be an injective $S$--module. The complex
  $\Hom{\crs{L}}{\s{I}}$ consists of injective $S$--modules, and it is
  exact as $\s{I}\in\I[S] \subseteq \I$.  Let $\s{J}$ be any injective
  $S$--module; swap gives $\Hom[S]{\s{J}}{\Hom{\crs{L}}{\s{I}}} \is
  \Hom{\crs{L}}{\Hom[S]{\s{J}}{\s{I}}}$.  This complex is exact by
  \prpcite[(4.1.3)]{LWC}, as $\Hom[S]{\s{J}}{\s{I}} \eq
  \DHom[S]{\s{J}}{\s{I}} \in \F$ by \lemref{stabcpxnoeth}.
  
  \HGF: Let $\crs{L}$ be a complete resolution by finite free
  $R$--modules and $\s{F}$ be a flat $S$--module. The complex
  $\Hom{\crs{L}}{\s{F}}$ consists of flat $S$--modules and is exact by
  \prpcite[(4.1.3)]{LWC} as $\s{F}\in\F[S] \subseteq \F$.  Let $\s{J}$
  be any injective $S$--module; tensor evaluation \respartref{ev}{c}
  gives an \iso $\tp[S]{\s{J}}{\Hom{\crs{L}}{\s{F}}} \is
  \Hom{\crs{L}}{\tp[S]{\s{F}}{\s{J}}}$.  This complex is exact, as the
  module $\tp[S]{\s{F}}{\s{J}} \eq \Dtp[S]{\s{F}}{\s{J}} \in \I$ by
  \lemref{stabcpxnoeth}.
  
  \HGP: Let $\crs{L}$ be a complete resolution by finite free
  $R$--modules. If $\s{P}$ is a projective $S$--module, the complex
  $\Hom{\crs{L}}{\s{P}}$ consists of projective $S$--modules and is
  exact by \prpcite[(4.1.3)]{LWC}.  Let $\s{Q}$ be any projective
  $S$--module; Hom~evaluation \respartref{ev}{a} gives an \iso
  $\Hom[S]{\Hom{\crs{L}}{\s{P}}}{\s{Q}} \is
  \tp{\crs{L}}{\Hom[S]{\s{P}}{\s{Q}}}$.  By \lemref{stabcpxnoeth} we
  get $\Hom[S]{\s{P}}{\s{Q}} \eq \DHom[S]{\s{P}}{\s{Q}} \in \F$, and
  thus $\tp{\crs{L}}{\Hom[S]{\s{P}}{\s{Q}}}$ is exact by
  \prpcite[(4.1.3)]{LWC}.
  
  \HFB: First we assume (1). Let $\s{F}$ be flat over $S$, let $B$ be
  \gor injective over $R$ and consider \mbox{$\Hom{\s{F}}{B} \is
    \Hom[S]{\s{F}}{\Hom{S}{B}}$}. Since \mbox{$\P[S] \subseteq \P$}
  the module $\Hom{S}{B}$ is \gor injective over $S$ by \HPBB, and
  since \mbox{$\F[S]=\P[S]$}, it follows by \pgref{table1}\HFBB\ 
  applied to $\f=1_R$ that $\Hom[S]{\s{F}}{\Hom{S}{B}}$ is \gor
  injective over $S$.
  
  Next we assume (2). Let $\crs{I}$ be a complete injective resolution
  over $R$, then the complex $\Hom{\s{F}}{\crs{I}}$ consists of
  injective $S$--modules, and it is exact as $\s{F}$ belongs to $\F[S]
  \subseteq\P$.  For any injective $S$--module $\s{J}$ the complex
  $\Hom[S]{\s{J}}{\Hom{\s{F}}{\crs{I}}} \is
  \Hom{\tp[S]{\s{J}}{\s{F}}}{\crs{I}}$ is exact, as
  $\tp[S]{\s{J}}{\s{F}}\in\I$ by \resref{stabcpx}, if $S$ is
  noetherian, or by \lemref{stabcpxnoeth}, if $R$ is noetherian.
\end{prf}

%%% ACKNOWLEDGMENTS

\section*{Acknowledgments}
\label{sec:ack}
\noindent
We thank Srikanth Iyengar for letting us include \thmref{ABE} in this
paper and Hans-Bj{\o}rn Foxby for showing us \pgref{table2}\TPX.  It
is also a pleasure to thank the referee for thorough and useful
comments.

\enlargethispage*{0.5cm}
%%% BIBLIOGRAPHY

\bibliographystyle{amsplain}

\providecommand{\bysame}{\leavevmode\hbox to3em{\hrulefill}\thinspace}
\renewcommand{\MR}{\relax\ifhmode\unskip\space\fi MR }
                                %\providecommand{\MRhref}[2]{%
                                %  \href{http://www.ams.org/mathscinet-getitem?mr=#1}{#2}
                                %}
                                %\providecommand{\href}[2]{#2}

\end{document}